\documentclass[12pt,leqno]{amsart}
\usepackage{amsmath,stmaryrd}
\usepackage{amssymb}
\usepackage{amsthm}
\usepackage{epsf}
\usepackage{graphicx}
\usepackage{esint} 
\usepackage{dsfont}
\usepackage[pagebackref]{hyperref} 
\hypersetup{pdfpagemode=FullScreen,  colorlinks=true} 
\usepackage{enumerate} 

\numberwithin{equation}{section}

\newcommand{\ms}{\medskip}

\newcommand{\R}{\mathbb{R}}

\newcommand{\bN}{\mathbb{N}}

\newcommand{\dist}{\,\mathrm{dist}}

\newcommand{\sm}{\setminus}
\newcommand{\supp}{\mathrm{supp}}
\newcommand{\diam}{\mathrm{diam}}

\newcommand{\wt}{\widetilde}

\newcommand{\cW}{{\mathcal  W}}

\newcommand{\cF}{{\mathcal  F}}

\newcommand{\1}{{\mathds 1}}

\newcommand{\dr}{\partial}

\DeclareMathOperator{\diver}{div}

\theoremstyle{plain}
\newtheorem{theorem}[equation]{Theorem}
\newtheorem{lemma}[equation]{Lemma}

\newtheorem{proposition}[equation]{Proposition}

\newtheorem{definition}[equation]{Definition}

\theoremstyle{definition}

\theoremstyle{remark}
\newtheorem{remark}[equation]{Remark}

\newcommand{\re}{\mathbb{R}}

\newcommand{\bp}{\noindent {\it Proof}.\,\,}
\newcommand{\ep}{\hfill$\Box$ \vskip 0.08in}

\def\div{\mathop{\operatorname{div}}}

\newcommand{\Rd}{\color{red}}

\begin{document}

\title[Absolute continuity of the harmonic measure on low dim. rectifiable sets]{Absolute continuity of the harmonic measure on low dimensional rectifiable sets}

\author[Feneuil]{Joseph Feneuil}
\address{Joseph Feneuil. Dipartimento di Matematica, Universit\`a di Pisa, Largo Bruno Pontecorvo, 7, I-78349 Pisa, Italy}
\email{joseph.feneuil@dm.unipi.it}

\maketitle

\begin{abstract}
In the past decades, we learnt that uniform rectifiability is often a right candidate to go past Lipschitz boundaries in boundary value problems. If $\Omega$ is an open domain in $\R^n$ with mild topological conditions, we can even characterize the $n-1$ dimensional uniformly rectifiability of the boundary $\partial \Omega$ by the $A_\infty$-absolute continuity of the harmonic measure on $\partial \Omega$ with respect to the surface measure. 

In low dimension ($d<n-1$), David and Mayboroda tackled one direction of the above characterization, i.e. proved that if $\Gamma$ is a $d$-dimensional uniformly rectifiable set, then the harmonic measure (associated to an suitable degenerate elliptic operator) on $\Gamma$ is $A_\infty$-absolutely continuous with respect to the $d$-dimensional Hausdorff measure.

In the present article, we use a completely new approach to give an alternative and significantly shorter proof of David and Mayboroda's result. 
\end{abstract}

\ms\noindent{\bf Key words.}
Uniform rectifiability, $A_\infty$-absolute continuity, harmonic measure, low dimensional boundaries, degenerate elliptic operators.

\ms\noindent

\tableofcontents

\section{Introduction}

\subsection{History and Motivation}

The past decades have seen a considerable achievements at the intersection of harmonic analysis, PDEs, and geometric measure theory. The general idea of the research is to establish links between the geometry of the boundary of a domain $\Omega$ and the regularity of the solutions (or the well-posedness of boundary value problems). Let us give an example. Given an open domain $\Omega \subset \R^n$ and a function $g$ on the boundary $\partial \Omega$, the Dirichlet problem consist to find solutions (possibly in some weak sense and appropriate spaces) to $-\Delta u = 0$ in $\Omega$ and $u=g$ on $\partial \Omega$. If $\Omega$ is a $C^{k,\alpha}$-domain and $g\in C^{k,\alpha}(\overline{\Omega})$ for a $k\geq 2$, then it is well known that the solution $u$ to the aforementioned Dirichlet problem is also $C^{k,\alpha}(\overline{\Omega})$ (see e.g. \cite[Theorem 6.19]{GT}).

Of course, mathematicians have studied the Dirichlet conditions under weaker conditions on $\partial \Omega$ and $g$. The important discovery for us is the following equivalence: the solvability of the Dirichlet problem for all $g$ in $L^p$ and one large $p<+\infty$ is equivalent to the $A_\infty$-absolute continuity of the harmonic measure with respect to the surface measure, see e.g. \cite[Theorem 1.7.3]{KenigB}, the $A_\infty$-absolute continuity is a scale invariant quantitative version of the mutual absolute continuity. Using the harmonic measure instead of boundary value problem is convenient because captures the diffusion property of the Laplacian, but places more emphasis on the boundary $\dr \Omega$ than on the domain $\Omega$.

The first result that links the harmonic measure and the boundary is now more than a century old. In 1916, F. and M. Riesz showed that for simply connected domains in the complex plane with rectifiable boundary, the harmonic measure is absolutely continuous with respect to the arc length (see \cite{RR}). In 1936, M. Lavrent'ev established a scale invariant version of the Riesz brothers result (see \cite{Lv}). C. Bishop and P. Jones obtained a local version of the result in 1990 (\cite{BJ}), and showed that topological conditions are needed to ensure that the harmonic measure on a rectifiable boundary is absolutely continuous with respect to the arc length.

The result was also studied in higher dimensional spaces, namely $\R^n$ for $n\geq 3$. B. Dahlberg proved in 1977 that, for domains $\Omega$ with Lipschitz boundary, the harmonic measure is indeed absolutely continuous continuous with respect to the surface measure $\mathcal H^{n-1}|_{\dr \Omega}$ (see \cite{Da}). The topic underwent a lot of improvements in the next three decades, leading to finer and finer necessary and sufficient conditions, see for instance \cite{JK}, \cite{DJ}, \cite{Se}, \cite{Ba}, \cite{Wu}, \cite{Z}. The authors of \cite{HM1} (see also \cite{AHMNT}) obtained that under some conditions of topology, uniform rectifiability of the boundary implies that the harmonic measure is $A_\infty$-absolutely continuous with respect to the Hausdorff measure was also observed in \cite{HMU} (with topological assumptions) and \cite{AHM3TV}--\cite{HLMN} (without topological assumptions) that rectifiability is necessary to get absolute continuity of the harmonic measure. 
Several recent works (e.g. \cite{Azzam}, \cite{AHMMT})  looked then for the optimal topological condition to have in order to ensure, together with uniform rectifiability, the $A_\infty$-property  of the harmonic measure (or slightly weaker versions).
A lot of articles could be related to the present discussion, for instance a lot of work has been done to see to which extend we can replace - in the previous results - the harmonic measure by elliptic measures associated to elliptic operators other than the Laplacian (which is relevant since at least one approach to this problem is to study perturbations of the Laplacian for simpler domains). We apologize to all the other mathematicians that could have been cited here, we refer to \cite{DFMprelim2} for a more detailed presentation of the state of the art, and we send the non-specialist interested reader to \cite{ToroICM} and \cite{PipherICM} for nice presentations of the problems related to this topic.

\medskip

Guy David and Svitlana Mayboroda had the following thoughts. Would it be a way, for sets $\Gamma \subset \R^n$ of dimension $d<n-1$, to obtain a similar criterion of uniform rectifiability using harmonic measure? A positive answer would be a huge discovery, because most of the criterions of uniform rectifiability, in particular those pertaining to PDE, are limited to certain dimensions or codimensions. The question at that time started as a challenge, since the harmonic measure is a tool that - roughly speaking - can only ``see'' the parts of $\Gamma$ sets of dimension $d$ such that $n-2<d<n$.
The idea was thus to find a way to construct a probability measure by way of elliptic PDEs, that will play the role of the harmonic measure, for sets with low dimension. In for instance \cite{LN}, the authors used a non-linear $p$-Laplacian operator to solve this issue, but their goal was different from David and Mayboroda's objective. David and Mayboroda's approach was to use linear but degenerate elliptic operators $L$ in $\Omega := \R^n \setminus \Gamma$, that satisfy elliptic and boundedness conditions relative to a weight $w(x) = \dist(x,\Gamma)^{d+1-n}$ that takes the dimension $d$ of $\Gamma$ and the distance to the boundary into account.
These ideas led to the memoir \cite{DFMprelim}, where we developed an elliptic theory associated to the degenerate operators that we wanted to use. In particular, we prove that when $L:=-\div [w(x)A(x) \nabla]$ is a degenerate elliptic operator and $A(x)$ satisfies the classical elliptic and boundedness conditions, weak solutions to $Lu=0$ in $\R^n \setminus \Gamma$ satisfy De Giorgi-Nash-Moser estimates inside the domain and at the boundary. We can then define a probability measure $\omega^X_{L}$ on $\Gamma$ associated to $L$, and this measure $\omega^X_{L}$ has desirable properties such as the doubling property, the degeneracy property, the change of pole property. We do not want to talk much about those properties; indeed, they are only used for the proof of Lemma \ref{Main13} below which will not be repeated here because it can be found in previous works.

In \cite{DFMAinfty}, we continued our project by aiming for Dahlberg's result for sets with a low dimension, that is if $\Gamma$ is the graph of a Lipschitz function from $\R^d$ to $\R^{n-d}$, then the ``harmonic measure'' $\omega^X_L$ is $A_\infty$-absolutely continuous with respect to the $d$-dimensional Hausdorff measure. The main difficulty here is the fact that - even in the classical case - not all elliptic operator with bounded coefficients satisfy that $\omega^X_L$ is $A_\infty$-absolutely continuous\footnote{Here and in the sequel, the $A_\infty$ absolute continuity will always be with respect to the $d$-dimensional Hausdorff measure on $\Gamma$.} (see \cite{CFK,MM}). We had to make a choice for $L:=L_\Gamma$, which is simple enough and systematically defined for all sets $\Gamma \subset \R^n$ of dimension $d$. The survey \cite{DFMKenig} presents what we succeeded to do by 2018. But what you need to know for the article - in particular the choice of $L_\Gamma$ - is given in the next subsection.

Guy David and Svitlana Mayboroda extended in \cite{DM} the above result to all  uniformly rectifiable sets. That is, for such sets, the harmonic measure $\omega^X_{L_\Gamma}$ is $A_\infty$ absolutely continuous for all uniformly rectifiable sets $\Gamma \subset \R^n$ of dimension $d<n-1$. Contrary to the case $d=n-1$, we do not need any assumption of topology on our domain $\Omega := \R^n \setminus \Gamma$, since they are automatically verified when $\Gamma$ has dimension $d<n-1$ (the fact that no extra topology condition is needed can be related to the fact that we are unlikely to touch the boundary when we travel between two points of $\Omega$). 

In the current article, we purpose a shorter and simpler proof of the main theorem in \cite{DM}. Our result is established by a completely different method that exploits in a crucial manner on the fact that $\Gamma$ has a dimension $d<n-1$. Our methods are simple in nature, for instance they do not rely on the so called Corona decomposition, saw-tooth domains, extrapolation arguments, ... 

\subsection{Main result} \label{SS1.2} In this subsection, we want to properly introduce all the tools that we need for our main result, and later state our main theorem. We shall first talk about the uniform rectifiability, and then turn to the presentation of the degenerate elliptic operator that will substitute the Laplacian and that will be used to construct our harmonic measure.

\medskip

Let $\Gamma \subset \R^n$ be a Ahlfors regular set of dimension $d$, which means that $\Gamma$ is a closed set and there exists a measure $\sigma$ supported on $\Gamma$ and a constant $C_{\sigma} \geq 1$ such that 
\begin{equation} \label{defAR}
C_{\sigma}^{-1} r^{d} \leq \sigma (B(x,r)) \leq C_\sigma r^d
\end{equation}
for all $x\in \Gamma$ and $r>0$. The Ahlfors regularity is a property on the set $\Gamma$ rather than on the measure $\sigma$. Indeed, if a $\sigma$ satisfying \eqref{defAR} exists, then \eqref{defAR} is also satisfied when $\sigma$ is the $d$-dimensional Hausdorff measure on $\Gamma$, possibly with a larger constant $C_\sigma$.

The geometric assumption on $\Gamma$ in this article is uniformly rectifiability. Uniformly rectifiable sets were introduced by David and Semmes, and equivalent definitions are given in \cite{DS1,DS2}. The characterization of uniform rectifiability that would be closer to the one of rectifiability is probably the one stating 
\begin{equation} \label{defBPLI}
\text{$\Gamma$ is uniformly rectifiable if $\Gamma$ has big pieces of Lipschitz images,}
\end{equation}
that is, if $\Gamma$ is Ahlfors regular \eqref{defAR}, and there exists $\theta, M>0$ such that, for each $x\in \Gamma$ and $r>0$, there is a Lipschitz mapping $\rho$ from the ball $B(0,r) \subset \R^d$ into $\R^n$ such that $\rho$ has Lipschitz norm $\leq M$ and 
\[\sigma(\Gamma \cap B(x,r) \cap \rho(B_{\R^d}(0,r))) \geq \theta r^d.\]
 In this article, we shall only rely on the characterization of uniform rectifiable sets by Tolsa $\alpha$-numbers that we present now.  
We denote by $\Xi$ the set of affine $d$-dimensional planes in $\R^n$. Each plane $P\in \Xi$ is associated with a measure $\mu_P$, which is the restriction to $d$-dimensional Hausdorff measure to $P$ (i.e. the Lebesgue measure on the plane). A {\bf flat measure} is a measure $\mu$ that can be written $\mu = c\mu_P$ where $c$ is a positive constant and $P\in \Xi$. The set of flat measures is called $\mathcal F$. 
We need Wasserstein distances to quantify the difference between two measures, and we shall use it to measure how far a measure $\sigma$ is from a flat measure.

\begin{definition} \label{D6.1}
For $x\in \R^n$ and $r > 0$, denote by $Lip(x,r)$ the set of $1$-Lipschitz functions $f$ supported in $\overline{B(x,r)}$, that is the set of functions $f : \R^n \to \R$ such that $f(y)=0$ for $y\in \R^n \sm B(x,r)$ and $|f(y)-f(z)|\leq |y-z|$ for $y,z\in \R^n$. The normalized Wasserstein distance between two measures $\sigma$ and $\mu$ is then
\begin{equation} \label{a6.2}
\dist_{x,r}(\mu,\sigma) = r^{-d-1} \sup_{f\in Lip(x,r)} \Big|\int f d\sigma - \int f d\mu\Big|.
\end{equation}
The distance to flat measures is defined by
\begin{equation} \label{a6.3}
\alpha_\sigma(x,r) = \inf_{\mu \in \cF}\dist_{x,r}(\mu,\sigma).
\end{equation}
\end{definition}

Observe that $\alpha_\sigma$ is uniformly bounded, i.e. there exists a constant that depends only on $d$, $n$, and $C_\sigma$ such that for all $x\in \Gamma$ and $r>0$, 
\begin{equation} \label{alphabdd}
\alpha_\sigma(x,r) \leq C_\sigma.
\end{equation}
The result above is quite classical. Take a flat measure $\mu$ supported outside $B(x,r)$, and we can see that $\alpha_\sigma(x,r) \leq r^{-d-1} \sup_{f\in Lip(x,r)} \Big|\int f d\sigma\Big|$. But since $f$ is 1-Lipschitz supported in $B(x,r)$, the function $f$ is bounded by $r$, which leads to the fact that $\alpha_\sigma(x,r) \leq r^{-d} \sigma(B(x,r)) \leq C_\sigma$ as desired.

Let $\Gamma$ be a $d$-Ahlfors regular set, and $\sigma$ be a measure that satisfies \eqref{defAR}. If $\Gamma$ is uniformly rectifiable, then there exists a constant $C_0>0$ that depends only on $\sigma$ such that 
\begin{equation} \label{defUR}
\int_{0}^r \int_{\Gamma \cap B(x,r)} |\alpha_\sigma(y,s)|^2 \, d\sigma(y) \, \frac{ds}{s} \leq C_0 \sigma(B(x,r)) \qquad \text{ for } x\in \Gamma \text{ and } r>0.
\end{equation}
The above statement is even a characterization of uniform rectifiability (see Theorem 1.2 in \cite{Tolsa09}). Tolsa's characterization of uniform rectifiability is given with dyadic cubes, but one can easily check that our bound \eqref{defUR} is equivalent to Tolsa's one.

\medskip

Now, we present the elliptic theory associated to our problem. Define $\Omega := \R^n \setminus \Gamma$, where $\Gamma$ is a $d$-Ahlfors regular set with $d < n-1$. The set $\Omega$ will serve as our domain in which we study elliptic equations. Because of the thin boundary, the domain $\Omega$ automatically satisfies the Harnack chain condition (quantitative connectedness), see Lemma 2.1 \cite{DFMprelim}. Moreover Lemma 11.6 in \cite{DFMprelim} entails the existence of a constant $C$ that depends only on $C_\sigma$ and $n-d>1$ such that for any $x\in \Gamma$ and $r>0$, we can find a point $A_{x,r}$ such that  
\begin{equation} \label{corkscrew}
C^{-1} r \leq \dist(A_{x,r},\Gamma) \leq |A_{x,r}-x| \leq Cr.
\end{equation}
So contrary to the case where $d=n-1$ (\cite{HM1}), we don't need to assume those topological hypotheses.

When $\Gamma$ is of codimension at least 2, a weak solution to $-\Delta u = 0$ in $\Omega:=\R^n \setminus \Gamma$ is also a weak solution to $-\Delta u = 0$ in $\R^n$. So in particular, we cannot impose any non-smooth data on $\Gamma$, and we cannot define the harmonic measure on $\Gamma$. 
In \cite{DFMprelim}--\cite{DFMprelim2}, the three authors developed an elliptic theory associated to a domain $\Gamma$ by considering degenerate  elliptic operators that takes the dimension of $\Gamma$ into account. In the present article, we consider the operators $L_{\beta,\gamma}$, $\beta >0$ and $\gamma \in (-1,1)$ defined as
\begin{equation} \label{defLbg}
L_{\beta,\gamma} := - \div (D_\beta)^{d+1+\gamma-n} \nabla.
\end{equation}
where $D_\beta$ is defined on $\Omega$ as
\begin{equation} \label{defDb}
D_\beta (X) := \left( \int_\Gamma |X-y|^{-d-\beta} d\sigma(y) \right)^{-1/\beta}.
\end{equation}
and $\sigma$ is the measure on $\Gamma$ introduced in \eqref{defAR}.
Lemma 5.1 in \cite{DFMAinfty} shows that, when $\Gamma$ is $d$-Ahlfors regular, 
\begin{equation} \label{Dbdist}
C^{-1} \dist(X,\Gamma) \leq D_\beta \leq C \dist(X,\Gamma) \qquad \text{ for } X\in \Omega,
\end{equation}
where $C>0$ depends only on $n$, $d$, $\beta$, and $C_\sigma$. In view of the above estimate, we can extend the definition of $D_\beta$ to all $\R^n$ by setting $D_\beta(x) = 0$ if $x\in \Gamma$. Moreover, it shows that the operator $L_{\beta,\gamma}$, for $\beta >0$ and $\gamma \in (-1,1)$ enters the scope of the theory written in \cite{DFMprelim2} (see the discussion in paragraph 3.3 from \cite{DFMprelim2} when $\gamma \neq 0$, see also \cite{DFMprelim} for the case $\gamma = 0$). 

For the rest of the article, we say that $u$ is a weak solution to $L_{\beta,\gamma} = 0$ if 
\begin{equation} \label{weaksol}
\int_\Omega (\nabla u \cdot \nabla \varphi) \, D_\beta^{d+1+\gamma-n} = 0 \qquad \text{ for } \varphi \in C^\infty_0(\Omega).
\end{equation}
In the integral above, we didn't specify that we integrate with respect to the $n$-dimensional Lebesgue measure. For the rest of the article, to lighten the notation, an integral without measure will always be an integral against the $n$-dimensional Lebesgue measure. The precise definition of the harmonic measure, as constructed in \cite{DFMprelim2} is then:

\begin{definition} \label{defhm}
For each $X\in \Omega$, we can define a unique probability measure $\omega^X:=\omega^X_{\beta,\gamma}$ on $\Gamma$ with the following properties. For any compactly supported continuous function $g$ on $\Gamma$, the function $u_g$ defined as
\[ u_g(X) = \int_\Gamma g(y) d\omega^X(y)\]
is a weak solution to $L_{\beta,\gamma}u=0$ in $\Omega:=\R^n \setminus \Gamma$, which in addition is continuous on $\R^n$ and is equal to $g$ on $\Gamma$.
\end{definition}

The goal of the article is to obtain the following result.

\begin{theorem} \label{Main1}
Let $\Gamma \subset \R^n$ be a $d$-Ahlfors regular uniformly rectifiable set with $d < n-1$, and let $\sigma$ be an Ahlfors regular measure on $\Gamma$ that satisfies \eqref{defAR}. Define $L_{\beta,\gamma}$ as in \eqref{defLbg}. 
Then the associated harmonic measure satisfies $\omega_{\beta,\gamma}^X \in A_\infty(\sigma)$.
This means that for every choice of $\epsilon \in (0,1)$, 
there exists $\delta \in (0,1)$, that depends only on $C_\sigma$, $C_0$, $\epsilon$, $n$, $d$, $\beta$, and $\gamma$,
such that for each choice of $x\in \Gamma$, $r>0$, a Borel set $E\subset B(x,r) \cap \Gamma$, and a corkscrew point $X = A_{x,r}$ as in \eqref{corkscrew}, 
\begin{equation} \label{AiTh}
\frac{\omega^X_{\beta,\gamma}(E)}{\omega^X_{\beta,\gamma}(B(x,r) \cap \Gamma)} < \delta \Rightarrow \frac{\sigma(E)}{\sigma(B(x,r) \cap \Gamma)} < \epsilon.
\end{equation}
\end{theorem} 

Observe that \eqref{AiTh} implies that $\omega^X_{\beta,\gamma}$ is absolutely continuous with respect to $\sigma$ and, as raised earlier in the introduction, the $A_\infty$ property can be seen as a quantitative scale invariant version of the absolute continuity. 

The fact that two measures $\mu$, $\nu$ satisfy $\mu \in A_\infty(\nu)$ has several characterizations, as mentioned in \cite[Theorem 1.4.13]{KenigB}. 
It is worth mentioning that, contrary to what the notation suggests, the $A_\infty$ property is an equivalence relationship, that is $\mu \in A_\infty(\nu)$ is actually the same as $\nu \in A_\infty(\mu)$. 
Moreover, $\mu \in A_\infty(\nu)$ is equivalent to reverse H\"older bounds at all scales on the kernel $\frac{d\mu}{d\nu}$.

\subsection{Steps of the proof of Theorem \ref{Main1}}
We recall for a last time that $\Gamma \subset \R^n$ denotes a Ahlfors regular set of dimension $d\leq n-1$, and that $\Omega:= \R^n\setminus \Gamma$ is its complement. Moreover, $\sigma$ stands for an Ahlfors regular measure that satisfies \eqref{defAR}. This notation will stand for the rest of the article.

The main result holds only when $d<n-1$, but some intermediate results pertaining to the geometry of uniformly rectifiable sets will be true for any integer $d < n$.

In addition, we shall use the convenient symbols $\lesssim$ and $\eqsim$. The inequality $A \lesssim B$ means that $A$ is smaller than a constant times $B$, with a constant that depends on parameters that are either recalled or obvious from context. Similarly, $A\eqsim B$ means that $A\lesssim B$ and $B \lesssim A$.

\medskip

We shall prove the $A_\infty$-property of the harmonic measure via the following result.

\begin{lemma} \label{Main13}
Let $d< n-1$ and take $\gamma \in (-1,1)$. Consider the operator $L:= -\div D_\beta^{d+1+\gamma-n} \mathcal A \nabla $, where $\beta >0$ and $\mathcal A$ is a (measurable real) matrix function on $\Omega$ that satisfies the usual elliptic conditions, that is
\begin{enumerate}[(i)]
\item $|\mathcal A(X)\xi\cdot \zeta| \leq C_2 |\xi||\zeta|$,
\item $|\mathcal A(X)\xi \cdot \xi| \geq (C_2)^{-1} |\xi|^2$.
\end{enumerate}
Assume that we can find $K$ such that for any ball $B \subset \R^n$ centered on $\Gamma$ and any Borel set $H \subset \Gamma$, the solution $u_H$ defined by $u_H(X):=\omega^X_L(H)$ satisfies 
 \begin{equation} \label{main13a}
 \int_{B} |\nabla u_H|^2 D_\beta^{d+2-n} \leq K \sigma(B).
\end{equation}
Then the harmonic measure $\omega^X_L$ is $A_\infty(\sigma)$ in the sense given in Theorem \ref{Main1}.
\end{lemma}

In the last lemma, the choices of $\beta$ and $\mathcal A$ does not really matter. The conditions we gave are the ones that allow us to fall in the scope of \cite{DFMprelim2} and thus that ensure the existence and the properties of the harmonic measure, namely the non-degeneracy of the harmonic measure, the fact that $\omega^X_{\beta,\gamma}$ is a doubling measure, and the change of pole property (these results are respectively Lemma 15.1, Lemma 15.43, and Lemma 15.61 in \cite{DFMprelim2}). Furthermore, an analogue of Lemma \ref{Main13} exists also when $d=n-1$; and in this case, $\Gamma$ is the boundary of an open domain $\Omega$ satisfying extra topological conditions (Harnack chains and corkscrew points) which are automatically true when $d<n-1$. But since we do not what to give details for a situation that we do not need, we excluded the case $d=n-1$. 

The demonstration of lemma \ref{Main13} will not be given here. Even if the lemma is stated in a slightly different context than what you can currently find in the literature, this type of result is not surprising - it became classical in the past years for experts in the field and known as ``the BMO-solvability implies $A_\infty$-absolute continuity of the harmonic measure'' - and the proof would just a small variation of what have already been done. Theorem 8.9 in \cite{DFMAinfty} which deals with the case where $\gamma = 0$ and $\Omega = \R^n \setminus \R^d$ (but those conditions are not relevant for the proof) and is stated in a similar manner as our lemma. Earlier references are Theorem 3.2 in \cite{KKiPT} and Theorem 1.3 in \cite{DPP2015}. See also Theorem 4.22 \cite{MZ}, and in \cite{HLM}, the proof of Lemma 5.24 and how it is paired with Theorem 1.3 to prove Theorem 5.30. 

\medskip

Our objective switched now to the proof of \eqref{main13a}, which says that $|\nabla u_H|^2 D_\beta^{d+2-n} dX$ is a Carleson measure. With the appearance of Carleson measure, let us introduce the following condition, that will play a crucial role in the article.

\begin{definition} \label{defCM}
The function $f$ on $\Omega$ satisfies the Carleson measure condition if $f\in L^\infty(\Omega)$, and the quantity $|f(X)|^2 \dist(X,\Gamma)^{d-n} dX$ is a Carleson measure, that is if for any $x\in \Gamma$ and $r>0$, there holds 
\begin{equation} \label{defCMa}
\int_{B(x,r)} |f(X)|^2 \dist(X,\Gamma)^{d-n}dX \leq C \sigma(B(x,r))
\end{equation}
with a constant $C>0$ independent of $x$ and $r$.

For short, we shall write that $f \in CM$ or $f\in CM(C)$ when we want to refer to the constant in \eqref{defCMa}.
\end{definition}

Recall that $\dist(X,\Gamma) \eqsim D_\beta$, where $D_\beta$ is the quantity defined in \eqref{defDb}. We will use and abuse of this equivalence in all our article. In particular, we shall use or prove the Carleson measure condition with the quantity $D_\beta^{d-n}(X)$ instead of $\dist(X,\Gamma)^{d-n}$, and $\beta$ will be chosen to fit our purpose.

In order to prove \eqref{main13a}, we shall use Carleson perturbations in the spirit of Kenig and Piper in \cite{KePiDrift}, as we already did in \cite{DFMAinfty}.

\begin{lemma} \label{Main12}
Let $d \leq n-1$ and take $\gamma \in (-1,1)$. Consider the operator $L:= -\div D_\beta^{d+1+\gamma-n} \mathcal A \nabla $, where $\beta >0$ and $\mathcal A$ is a matrix function on $\Omega$ that satisfies the usual elliptic conditions given in Lemma \ref{Main13}. 
Assume that we can find a scalar function $b$ and a vector function $\mathcal V$, both defined on $\Omega$, such that 
\begin{equation} \label{main12z}
\div [(b\mathcal A^T \nabla D_\beta + \mathcal V) D_\beta^{d+1-n}] = 0
\end{equation}
and such that $b$ and $\mathcal V$ satisfy
\begin{enumerate}[(i)]
\item $C_1^{-1} \leq b \leq C_1$,
\item $D_\beta \nabla b \in CM(C_1)$,
\item $|\mathcal V| \leq C_1$,
\item $\mathcal V \in CM(C_1)$,
 \end{enumerate}
for some constant $C_1>0$.

Then for any ball $B\subset \R^n$ centered on $\Gamma$ and any weak solution $u$ to $Lu = 0$ in $2B$,  one has
 \begin{equation} \label{main12a}
 \int_{B} |\nabla u|^2 D_\beta^{d+2-n} \leq C\left(\sup_{2B} |u|^2\right) \sigma(B),
 \end{equation}
 where $C$ depends only on $C_\sigma$, $C_1$, $\beta$, $\gamma$, $n$, and $d$. In particular \eqref{main13a} holds.
\end{lemma}

The expression \eqref{main12z} has to be taken in a weak sense, that is we shall use that for any test function $\varphi \in W^{1,1}(\Omega,\R^n)$ with compact support in $\Omega$, one has
\begin{equation} \label{main12y}
\int_\Omega \nabla \varphi \cdot (b\mathcal A \nabla D_\beta + \mathcal V) D_\beta^{d+1-n} = 0.
\end{equation}
Moreover, \eqref{main12a} is actually a typical $S<N$ estimate. For a ball $B$ centered on $\Gamma$, and $x\in \Gamma$, we define the cone in $B$ with vertex $x$ as $\gamma^B(x) := \{X\in B, \, |X-x| \leq 2\dist(X,\Gamma)\}$
and then, for a function $u$ defined on $\Gamma$
\[N^B(u)(x) = \sup_{\gamma^B(x)} |u|.\]
We actually prove the following stronger version of \eqref{main12a}: 
 \begin{equation} \label{main12b}
 \int_{B} |\nabla u|^2 D_\beta^{d+2-n} \leq C\|N^{2B}(u)\|_{L^2(6B)}^2.
 \end{equation}
 
Lemma \ref{Main12} is probably the key result, but at the same time, its proof is elementary and uses classical computations. Theorem 1.32 in \cite{DFMAinfty} states a similar result when $\Gamma = \R^d$. The proof in \cite{DFMAinfty} relies on the fact that $|t|$ is a solution to $L_0u=0$ with $L_0$ being a `Carleson perturbation' of the considered operator $L$.
In our case, the analogue of $|t|$ is $D_\beta$ and $L_0 u = \div [D_\beta^{d-n} (bD_\beta \mathcal A \nabla + \mathcal V) u]$.

Of course, the above lemma alone does not look very appealing. What is the point of considering the assumption \eqref{main12z}, which relies on the existence of two quantities $b$ and $\mathcal V$ that may be impossible to find? But its combination with the next geometrical result makes the magic happen.

\begin{lemma} \label{Main11a}
Let $\Gamma$ be uniformly rectifiable, i.e. \eqref{defUR} is verified, of dimension $d<n-1$. Let $\beta >0$.  Then there exist a scalar function $b$ and a vector function $\mathcal V$, both defined on $\Omega$, such that 
\[\int_{\Gamma} |X-y|^{-n}(X-y) d\sigma(y) = (b  \nabla D_\beta + \mathcal V) D_\beta^{d+1-n} \qquad \text{ for } X\in \Omega\]
\begin{enumerate}[(i)]
\item $C_1^{-1} \leq b \leq C_1$,
\item $D_\beta \nabla b \in CM(C_1)$,
\item $|\mathcal V| \leq C_1$,
\item $\mathcal V \in CM(C_1)$,
 \end{enumerate}
 where $C_1$ is a constant that depends only on $C_\sigma$, $C_0$, $\beta$, $n$, and $d$.
\end{lemma}

One extra observation is needed. The quantity $\int_{\Gamma} |X-y|^{-n}(X-y) d\sigma(y)$ is divergence free. Indeed, we can locally interchange derivative and integral by using the dominated convergence theorem in order to get
\begin{equation} \label{divH}
\div \int_{\Gamma} |X-y|^{-n}(X-y) d\sigma(y) = \int_\Gamma \diver_X [ |X-y|^{-n}(X-y) ] d\sigma(y) = 0.
\end{equation}
Therefore, Lemma \ref{Main11a} gives us exactly what we need to apply Lemma \ref{Main12}.

The limitation $d<n-1$ comes from the fact that the quantity $\int_{\Gamma} |X-y|^{-n}(X-y) d\sigma(y)$ is not defined when $d=n-1$. Nothing stops Lemma \ref{Main12} to be valid in the case $d=n-1$, but our lack of substitute for $\int_{\Gamma} |X-y|^{-n}(X-y) d\sigma(y)$ in co-dimension 1 is why we believe that our proof cannot be (easily) adapted to the classical co-dimension 1 case. 

\medskip

The next result is an interesting variant of Lemma \ref{Main11a}, which will be proved in Section \ref{SUR}, and used to prove Lemma \ref{Main11a}.

\begin{lemma} \label{Main11}
Let $\Gamma$ be uniformly rectifiable, i.e. \eqref{defUR} is verified, of dimension $d<n$. Then for any $\alpha,\beta >0$, the quantity $\dist(X,\Gamma) \nabla [D_{\beta}/D_{\alpha}]$ satisfies the Carleson measure condition with a constant that depends only on $C_\sigma$, $C_0$, $\alpha$, $\beta$, $n$, and $d$.
\end{lemma}


We conclude our section by saying that the reader is welcome to check that the combination of Lemmas \ref{Main11a}, \ref{Main12}, and \ref{Main13} easily implies Theorem \ref{Main1}. As a consequence, the rest of the article is solely devoted to the proofs of Lemma \ref{Main12} and Lemma \ref{Main11a} (in this order). The two demonstrations use different techniques, and the sections can be read independently. 

\bigskip

\noindent {\bf Acknowledgements:} The author would like to thank Svitlana Mayboroda and Guy David for many fruitful discussions and for providing an early presentation of Corollary 6.6 in \cite{DEM} from which the ideas of the present article are derived.

\section{Proof of Lemma \ref{Main12}}
\label{SP12}

In this section, $d$ can be any real value in $(0,n-1]$. The first step in the proof of Lemma \ref{Main12} is to establish a Carleson inequality. The proof Carleson inequality exists in $\R^{n+1}_+$ (see \cite{Stein93}), or in $\mathcal M \times (0,+\infty)$ where $\mathcal M$ is a manifold (see \cite{Russ}). We do not pretend that our context is more complicated or even that the arguments of the proof are different, but we could not pinpoint a good reference, so we believe that it was now a good opportunity to discuss (and sketch the proof) about the Carleson inequality on Ahlfors regular sets.

We start with the definition. We say that a function $f$ on $\Omega$ describes a Carleson measure if $f$ is a Borel measurable function and if $f(X)\dist(X,\Gamma)^{d-n}dX$ is a Carleson measure, that is if there exists $C>0$ such that for any ball $B$ centered on $\Gamma$,
\begin{equation} \label{defCM1}
\int_{B} |f(X)| \dist(X,\Gamma)^{d-n}dX \leq C \sigma(B).
\end{equation}
The quantity $\|f\|_{CM1}$ denotes the smallest constant that satisfies \eqref{defCM1} for any ball $B$ centered on $\Gamma$.
Then we need cones $\gamma(x)$ with vertex $x\in \Gamma$ defined by
\begin{equation} \label{defgamma}
\gamma(x) := \{X\in \Omega, \, |X-x| \leq 2 \dist(X,\Gamma)\};
\end{equation}
the constant $2$ in the definition of the cones $\gamma(x)$ does not matter. Any fixed constant $\alpha >1$ would do (and constants in the incoming estimates will then also depend on $\alpha$). The non-tangential maximal function $N$ is
\begin{equation} \label{N}
N(u)(x) := \sup_{\gamma(x)} |u|;
\end{equation}
if - say - $u$ is a continuous bounded function on $\Omega$ and $x\in \Gamma$.

We need the Carleson inequality.

\begin{proposition} \label{CarlesonI}
Let $f$ be a function on $\Omega$ such that $f(X) \dist(X,\Gamma)^{d-n} dX$ is a Carleson measure. There exists a constant $C>0$ that depends only on $C_\sigma$ such that for any continuous bounded function $u$ on $\Omega$,
\begin{equation} \label{CarlesonIa}
\left| \int_\Omega u(X) f(X) \dist(X,\Gamma)^{d-n} dX \right| \leq C \|f\|_{CM1} \int_\Gamma N(u) \, d\sigma.
\end{equation}
In particular, if $f\in CM(C_1)$ (see Definition \ref{defCM}), we easily deduce that 
\[\int_\Omega |u(X)|^2 |f(X)|^2 \dist(X,\Gamma)^{d-n} dX \leq C C_1 \|N(u)\|^2_{L^2(\Gamma,\sigma)}.\]
\end{proposition}

\bp The second part of the proposition is immediate from the first part. Without surprise the proof of \eqref{CarlesonIa} is the same (with obvious modifications) as the one in \cite[Section II.2.2, Theorem 2]{Stein93}, which treats the case $\Omega = \R^n_+$.
\ep

We have the proper tools to prove Lemma \ref{Main12}

\medskip

\noindent {\em Proof of Lemma \ref{Main12}.}
The result is a local one, so we use cut-off functions. Take $\psi \in C^\infty_0(\R)$ be such that $\psi \equiv 1$ on $[-1,1]$, $\psi$ is compactly supported in $(-2,2)$, $0\leq \psi \leq 1$, and $|\psi'|\leq 2$. Let $B=B(x,r)$ a ball centered on the boundary, and $\epsilon>0$. We define the function $\phi_{B,\epsilon}$ on $\Omega$ by
\begin{equation} \label{defphi}
\phi_{B,\epsilon}(X) :=  \underbrace{\psi\left(\frac{\dist(X,B)}{10\dist(X,\Gamma)}\right)}_{\psi_1(X)} \underbrace{\psi\left(\frac{2\dist(X,B)}{r}\right)}_{\psi_2(X)} \underbrace{\psi\left(\frac{\epsilon}{\dist(X,\Gamma)}\right)}_{\psi_3(X)} 
\end{equation}
Let us list few properties of $\phi_{B,\epsilon}$. We have 
\begin{equation} \label{phi1}
\phi_{B,\epsilon}(X) = 1 \quad \text{$X \in B$, $\dist(X,\Gamma) \geq \epsilon$}.
\end{equation}
In addition, the function  $\phi_{B,\epsilon}$ is  supported on 
\begin{equation} \label{suppphi}
\supp\,  \phi_{B,\epsilon} \subset \{X \in 2B, \, \dist(X,B) \leq 20\dist(X,\Gamma), \, \dist(X,\Gamma) \geq \epsilon/2\},
\end{equation}
At last, we want to bound its gradient. With the help of \eqref{suppphi}, we deduce that 
\[ |\nabla \phi_{B,\epsilon}| \leq \frac{10}{\dist(X,\Gamma)} \left( \1_{\supp \, \nabla \psi_1} + \1_{\supp \,  \nabla \psi_2} + \1_{\supp \, \nabla \psi_3} \right)  \]
We quickly observe that  
\[\begin{split}
\supp\, \nabla \psi_1 & \subset \{10\dist(X,\Gamma)\leq \dist(X,B) \leq 20 \dist(X,\Gamma)\} \\
\supp \, \nabla \psi_3 & \subset \{ \epsilon/2 \leq \dist(X,\Gamma) \leq \epsilon\} \\
\supp \, \nabla \psi_2 & \subset \{X\in 2B, \, r/2 \leq \dist(X,B)\}.
\end{split}\]
Together with the facts that $20 \dist(X,\Gamma) \geq \dist(X,B)$ and $\dist(X,\Gamma) \leq |X-x| \leq 2r$ when $X \in \supp \, \phi_{B,\epsilon}$, we obtain that
\begin{equation} \label{gradphi}
|\nabla \phi_{B,\epsilon}| \leq \frac{100}{\dist(X,\Gamma)} \left[ \1_{E_1} + \1_{E_2} + \1_{E_3}\right],
\end{equation}
where 
\begin{equation} \label{defE1}
E_1:= \{X\in 2B, \, 10\dist(X,\Gamma) \leq \dist(X,B) \leq 20 \dist(X,\Gamma)\},
\end{equation}
\begin{equation} \label{defE2}
E_2:= \{X\in 2B, \, r/40 \leq \dist(X,\Gamma) \leq 2r\},
\end{equation}
and
\begin{equation} \label{defE3}
E_3:= \{X\in 2B, \, \epsilon/2 \leq \dist(X,\Gamma) \leq \epsilon\}.
\end{equation}

We claim that $\1_{E_1}$, $\1_{E_2}$, and $\1_{E_3}$ satisfy all the Carleson measure condition, that is for any $(y,s)\in \Gamma \times (0,+\infty)$, 
\begin{equation} \label{claimEi}
\int_{B(y,s)} |\1_{E_1} + \1_{E_2} + \1_{E_3}|^2 \dist(X,\Gamma)^{d-n} dX \leq Cs^d.
\end{equation} 
Note that if we prove the claim, we also prove the same estimate without the square power.  We shall demonstrate the claim separately for each $E_i$. First, Fubini's identity and the Ahlfors regularity of $\sigma$ entail, for $i\in \{1,2,3\}$, that
\begin{equation} \label{claimEi1}
\int_{B(y,s)} |\1_{E_i}|^2 \dist(X,\Gamma)^{d-n} dX \leq C \int_{B(y,10s)} \left(\int_{\gamma(z)} \1_{E_i}(X) \dist(X,\Gamma)^{-n} dX\right) d\sigma(z).
\end{equation} 
where the cone $\gamma(x)$ is the one from \eqref{defgamma}. Therefore, it is enough to prove that for all $z\in \Gamma$, and for $i\in \{1,2,3\}$ one has 
\begin{equation} \label{claimEi2}
\int_{\gamma(z)} \1_{E_i}(X) \dist(X,\Gamma)^{-n} dX \lesssim 1.
\end{equation} 
On $E_3$, we have $\dist(X,\Gamma) \eqsim \epsilon$, so 
\[\int_{\gamma(z)} \1_{E_3}(X) \dist(X,\Gamma)^{-n} dX \lesssim \epsilon^{-n} \int_{\gamma(z) \cap E_3} dX \lesssim \epsilon^{-n} |B(z,2\epsilon)| \lesssim 1.\]
The estimate \eqref{claimEi2} for $i=3$ follows, so is the claim \eqref{claimEi} for $E_3$.
The claim \eqref{claimEi} for $E_2$ is similar to $E_3$, and is left to the reader. We turn to the proof of \eqref{claimEi2}, that implies \eqref{claimEi}, for $E_1$. Let $X \in \gamma(z) \cap E_1$. Having $X\in E_1$ means that 
\begin{equation} \label{distBGamma}
\dist(X,B) \leq 20\dist(X,\Gamma) \leq 2\dist(X,B).
\end{equation}
In addition, $X\in \gamma(z)$ means that $|X-z| \leq 2\dist(X,\Gamma)$. Combining this latter fact with the second inequality in \eqref{distBGamma} leads to 
\begin{equation} \label{distxz}
|X-z| \leq \frac15 \dist(X,B),
\end{equation}
and thus
\[ \frac45 \dist(X,B) \leq \dist(X,B) - |X-z| \leq \dist(z,B) \leq  |X-z| + \dist(X,B) \leq \frac65 \dist(X,B).\]
The bounds above and \eqref{distBGamma} allow us to compare $\dist(X,\Gamma)$ and $\dist(z,B)$. We have
\[ \frac1{24} \dist(z,B) \leq \frac1{20}\dist(X,B) \leq \dist(X,\Gamma) \leq \frac1{10} \dist(X,B) \leq \frac18 \dist(z,B)\]
which is very nice because we bounded $\dist(X,\Gamma)$ by quantities that do not depend on $X\in \gamma(z) \cap E_1$. All those estimates allow us to also say that $|X-z| \leq \frac14 \dist(z,B)$. As a consequence,
\[\int_{\gamma(z)} \1_{E_1}(X) \dist(X,\Gamma)^{-n} dX \lesssim \dist(z,B)^{-n} \int_{\gamma(z) \cap B(z,\frac14\dist(z,B))} dX \lesssim 1\]
The claims \eqref{claimEi2} and \eqref{claimEi} follow.

\medskip

Let us turn to the main part of the proof of Lemma \ref{Main12}. Let $u$ be a weak solution to $L_{\beta,\gamma} u = 0$ in $2B\cap \Omega$. We intend to prove that 
\begin{equation} \label{claimS<N} \begin{split}
\int_{\Omega} |\nabla u|^2 \phi_{B,\epsilon}^2 D_\beta^{d+2-n} \leq C & \int_\Gamma |N(u\1_{\supp\, \phi_{B,\epsilon}})|^2 \, d\sigma \\
& + C\left(\int_{\Omega} |\nabla u|^2 \phi_{B,\epsilon}^2 D_\beta^{d+2-n}\right)^\frac12 \left( \int_\Gamma |N(u\1_{\supp\, \phi_{B,\epsilon}})|^2 \, d\sigma \right)^\frac12,
\end{split} \end{equation}
with a constant $C>0$ that depends only on $d$, $n$, $C_\sigma$, and $C_1$.
Why is it enough? We used the cut-off function $\phi_{B,\epsilon}$, which is compactly supported in $\Omega$, and a weak solution $u$ to $Lu = 0$ \footnote{Therefore $u$ is in $W^{1,2}_{loc}(\Omega,dX) = W^{1,2}_{loc}(\Omega,D_\beta^{d+1-n}dX)$ by definition, and also in $L^\infty_{loc}(\Omega)$ because we have we have Moser estimates on solutions, see \cite{DFMprelim2}.}. 
Therefore all the quantities in \eqref{claimS<N} are finite. So the estimate \eqref{claimS<N} self-improves to
\[\int_{\Omega} |\nabla u|^2 \phi_{B,\epsilon}^2 D_\beta^{d+2-n} \leq C \int_\Gamma |N(u\1_{\supp\, \phi_{B,\epsilon}})|^2 \, d\sigma \]
The constant $C$ above is independent of $\epsilon$, so we can take the limit as $\epsilon \to 0$. The lemma follows then by the properties of $\phi_{B,\epsilon}$, in particular \eqref{phi1}. In order to prove the stronger bound \eqref{main12b}, we need to prove that $\|N(u\1_{\supp\, \phi_{B,\epsilon}})\|_{L^2} \leq \|N^{2B}(u)\|_{L^2(6B)}$, which is implied by the fact that $\gamma^{2B}(x):= \gamma(x) \cap 2B = \emptyset$ whenever $x\in \Gamma \setminus 6B$. 
But the latter is fairly straightforward. Indeed, if $X\in \gamma^{2B}(x)$, then \[\dist(x,B) \leq \dist(X,B) + |X-x| \leq r_B + 2\delta(X) \leq 5r_B,\]
where $r_B$ is the radius of the ball $B$ (centered on $\Gamma$). The bound \eqref{main12b} follows.

\medskip

It remains to establish \eqref{claimS<N}. For the rest of the proof, to lighten the notation, we write $\phi$ for $\phi_{B,\epsilon}$. Let $b,\mathcal V$ as in the assumptions of the Lemma. We  define $H_{n-d-1}$ as
\begin{equation} \label{HbNDb}
H_{n-d-1} := (b\mathcal A^T \nabla D_\beta + \mathcal V) D_{\beta}^{d+1-n},
\end{equation}
which is a quantity locally bounded in $\Omega$. The notation $H_{n-d-1}$ above may look a bit weird at this point (why not calling it simply $H$), but is consistent with the one in Section \ref{SUR}. We assume $\div[H_{n-d-1}] = 0$ in a weak sense, that is for any compactly supported test function $\varphi \in W^{1,1}(\Omega)$, we have
\begin{equation} \label{main12x}
\int_\Omega \nabla \varphi \cdot H_{n-d-1} = 0.
\end{equation}
Note also that $|H_{n-d-1}| \lesssim D_\beta^{d+1-n}$. The estimate is straightforward from the assumption on $b$ and $\mathcal V$ once you know that $|\nabla D_\beta| \lesssim 1$. The latter fact is not surprising, since $D_\beta$ is smooth and plays the role of a distance, and not very hard to prove from the definition either; but proof is postponed to \eqref{bdNDb}.

\medskip

Since $b \gtrsim 1$ and $\mathcal A$ is elliptic, we get
\[\int_{\Omega} |\nabla u|^2 \phi^2 D_\beta^{d+2-n}  \lesssim J:= \int_{\Omega}  (\mathcal A \nabla u \cdot \nabla u) \,  b\phi^2 D_\beta^{d+2-n}.\]
We use the product rule to force every term into the second gradient,
\[\begin{split}
J & = \int_{\Omega} \mathcal A \nabla u \cdot \nabla [u b\phi^2 D_\beta^{1-\gamma}] D_\beta^{d+1+\gamma-n} - 2 \int_{\Omega}   \mathcal A \nabla u \cdot \nabla \phi \,  bu \phi \,  D_\beta^{d+2-n} \\
& - \int_{\Omega}  \mathcal A \nabla u \cdot \nabla b \,  u \phi^2 \,  D_\beta^{d+2-n} - \int_{\Omega}  \mathcal A \nabla u \cdot \nabla [D_\beta^{1-\gamma}] \,  bu \phi^2 \,   D_\beta^{d+1+\gamma-n} \\
& := J_1 + J_2 + J_3 + J_4.
\end{split}\]
The term $J_1$ equals 0, because $u$ is a weak solution to $Lu = -\div D_\beta^{d+1+\gamma-n} \mathcal A \nabla u = 0$, and because \cite[Lemma 9.18]{DFMprelim} says that $v:= ub \phi^2 D_\beta^{1-\gamma}$ is a valid test function.
The terms $J_2$ and $J_3$ can be treated in a similar manner. Using the bounds on $b$ and $\mathcal A$, we have
\[\begin{split}
|J_2 + J_3| & \lesssim \int_\Omega |\nabla u| |u \1_{\supp \, \phi} \phi [|\nabla \phi| + |\nabla b|] D_\beta^{d+2-n} \\
& \lesssim \left(\int_\Omega |\nabla u|^2 \phi^2 D_\beta^{d+2-n} \right)^\frac12 \left( \int_\Omega [u\1_{\supp \, \phi}]^2  [|D_\beta\nabla \phi|^2 + |D_\beta \nabla b|^2] D_\beta^{d-n} \right)^\frac12\end{split}.\]
Recall that, thanks to \eqref{Dbdist}, the function $D_\beta$ can be used like $\dist(.,\Gamma)$, in particular in the Carleson inequality (Proposition \ref{CarlesonI}) that we intend to invoke. But first, let us verify that we have the Carleson measure condition we need. The fact that $D_\beta \nabla b \in CM$ is part of the assumption of the lemma. The fact that $D_\beta \nabla \phi \in CM$ is a consequence of \eqref{gradphi}, \eqref{claimEi}, and again \eqref{Dbdist}. Proposition \ref{CarlesonI} infers that
\[\int_\Omega [u\1_{\supp \, \phi}]^2 [|D_\beta\nabla \phi|^2 + |D_\beta \nabla b|^2] D_\beta^{d-n} \lesssim \int_\Gamma |N(u \1_{\supp \, \phi})|^2 \, d\sigma,\]
that is
\[|J_2 + J_3| \lesssim \left(\int_{\Omega} |\nabla u|^2 \phi^2 D_\alpha^{d+2-n}\right)^\frac12 \left( \int_\Gamma |N(u\1_{\supp\, \phi})|^2 \, d\sigma \right)^\frac12\]
which is perfect since the right-hand side above appears in the right-hand side of \eqref{claimS<N}.

The last term that we need to treat is $J_4$. We have
\[ J_4 = -(1-\gamma) \int_{\Omega}  \mathcal A \nabla u \cdot [b  D_\beta^{d+1-n}\nabla D_\beta] \,  u \phi^2 = (\gamma-1) \int_{\Omega}   \nabla u \cdot [  D_\beta^{d+1-n} b \mathcal A^T \nabla D_\beta] \,  u \phi^2 .\]
We use the relation \eqref{HbNDb} to get
\[\begin{split}
J_4 & = (1-\gamma) \int_{\Omega}  \nabla u \cdot \mathcal V \,  u \phi^2 \, D_{\beta}^{d+1-n} + (\gamma-1 ) \int_{\Omega}  \nabla u \cdot H_{n-d-1} \,  u \phi^2 \\
& := J_{41} + J_{42}.
\end{split}\]
The integral $J_{41}$ can be dealt with using the same computations as the ones we did for $J_{2} + J_3$, using the facts that $\mathcal V \in CM$. We are left with $J_{42}$, where we force all the terms in the first gradient with the product rule. We obtain that
\[\begin{split}
J_{42} & = \frac{\gamma-1}2 \int_\Omega \nabla[u^2\phi^2] \cdot H_{n-d-1} + (1-\gamma) \int_\Omega \nabla \phi \cdot H_{n-d-1} \, u^2 \phi \\
& := J_{421} + J_{422}.
\end{split}\]
The term $J_{421}$ is 0, due to \eqref{main12x}. As for $J_{422}$, similarly to what we did for $J_2$, we use \eqref{gradphi} and the fact that $|H_{n-d-1}| \lesssim D_\beta^{d+1-n}$ to get
\[\begin{split} J_{422} & \lesssim \int_\Omega [\1_{E_1} + \1_{E_2} + \1_{E_3}] \, [u^2 \1_{\supp \, \phi}] \, D_\beta^{d-n} \\
& \lesssim \int_\Omega |\1_{E_1} + \1_{E_2} + \1_{E_3}|^2 \, [u^2 \1_{\supp \, \phi}] \, D_\beta^{d-n}.\end{split}\]
 Proposition \ref{CarlesonI} and  \eqref{claimEi} prove that 
\[\begin{split}
J_{422} & \lesssim \|\1_{E_1} + \1_{E_2} + \1_{E_3}\|_{CM1} \int_{\Gamma} N(u^2\1_{\supp\, \phi}) \, d\sigma \\
& \lesssim \int_{\Gamma} |N(u\1_{\supp\, \phi})|^2 \, d\sigma
\end{split}\]
since $\|\1_{E_1} + \1_{E_2} + \1_{E_3}\|_{CM1}$ is bounded by a constant depending only on $d$, $n$, $C_\sigma$. The claim \eqref{claimS<N} and then the lemma follows. 
\ep

\section{Proof of Lemma \ref{Main11a}}

\label{SUR}

In this section $\Gamma$ is a $d$-Ahlfors regular set of dimension $d<n$\footnote{We prove a more general result (Lemma \ref{Main11b}) and our limitation $d<n-1$ in Lemma \ref{Main11a} only comes from the fact that Lemma \ref{Main11a} is Lemma \ref{Main11b} with $\alpha := n-d-1>0$.}. We deduce that $\Gamma$ is a closed non-empty set. So we can use a family of Whitney cubes $\cW$ as constructed in \cite{Stein93}. The diameter of $Q$ is written $\ell(Q)$, and we notice that the side length of $Q$ is $\ell(Q)/\sqrt n$. 

We record a few of the properties of $\cW$ that we shall need. The collection $\cW$ is the family of maximal dyadic cubes such that $20Q \subset \Omega$, that is we have
\begin{equation} \label{prW1}
20Q \subset \Omega \quad \text{ but } \quad 60Q \cap \Gamma \neq \emptyset.
\end{equation}
If $Q,R \in \cW$ are such that $2Q\cap 2R \neq \emptyset$, then $\ell(R) \in \{\ell(Q)/2,\ell(Q),2\ell(Q)\}$. Thus, $R$ is a dyadic cube that satisfies $R \subset 8Q$ and $\ell(R) \eqsim \ell(Q)$, and there are only a (uniformly) finite number of such cubes. This proves that there is a constant $K:=K(n)$ such that 
\begin{equation} \label{prW3}
\text{ the number of cubes $R\in \cW$ such that $2R \cap 2Q \neq \emptyset$ is at most $K$.}
\end{equation}

For each $Q\in \cW$, we pick once for all the article a point 
\begin{equation} \label{defxiQ}
\xi_Q \in 60Q \cap \Gamma.
\end{equation}
We write $B_Q$ for the ball $B_Q:= B(\xi_Q,2^5\ell(Q)) \supset Q$
and, using \eqref{a6.2}, we define the Wasserstein distance between two measures $\mu_1$ and $\mu_2$ relatively to $Q$ as
\begin{equation} \label{defBQ}
\dist_Q(\mu_1,\mu_2):= \dist_{\xi_Q,2^{10}\ell(Q)} (\mu_1,\mu_2).
\end{equation}
Consider
\begin{equation} \label{defaQ}
\alpha_Q:= \alpha_\sigma(\xi_Q,2^{10}\ell(Q)),
\end{equation}
and we have the following result.

\begin{lemma} \label{lemmuQ}
Let $\mu_Q:=c_Q\mu_{P_Q}$ be a flat measure that satisfies $\dist_{Q}(\mu_Q,\sigma) \leq 2 \alpha_Q$. 

There exists a small constant $\epsilon:=\epsilon(C_\sigma,d) >0$ such that if $\alpha_Q \leq \epsilon$, we have $\dist(\xi_Q,P_Q) \leq 5\ell(Q)$, $\dist(2Q,P_Q) \geq 5\ell(Q)/\sqrt n$, and $\epsilon \leq c_Q \leq 1/\epsilon$ .
\end{lemma}

\bp We shall prove the result by contraposition. Assume first that $\dist(\xi_Q,P_Q) \geq 5\ell(Q)$. In this case, we choose the 1-Lipschitz function 
\[\tilde f := \max\{5\ell(Q) - |X-\xi_Q|,0\}\]
in the definitions \eqref{a6.2} and\eqref{a6.3} to get
\[\left| \int \tilde f \, d\sigma - \int \tilde f \, d\mu_Q \right| \leq 2 [2^{10}\ell(Q)]^{d+1} \alpha_Q.\]
The function $\tilde f$ is always 0 on $P_Q$ and is at least $4\ell(Q)$ on $B(\xi_Q,\ell(Q))$. So the estimate above becomes 
\[ 4\ell(Q) \sigma(B(\xi_Q,\ell(Q))) \leq 2 [2^{10}\ell(Q)]^{d+1} \alpha_Q.\]
We use \eqref{defAR} to obtain a uniform lower bound on $\alpha_Q$. 

Assume now that $c_Q$ is smaller than a constant $\epsilon_0$ that depends only on $d$ and $C_\sigma$ and that will be chosen later. In this case we choose the 1-Lipschitz function $\hat f(X) = \max\{\ell(Q) - \dist(X, B_Q),0\}$ in \eqref{a6.2} to deduce that 
\begin{equation} \label{muQb}
\begin{split}
\alpha_Q & \gtrsim \ell(Q)^{-d-1} \left|\int \hat f \, d\sigma - \int \hat f \, d\mu\right| \\
& \gtrsim \ell(Q)^{-d} \left( \sigma(B_Q) - c_Q \mu_{P_Q}(2B_Q) \right).
\end{split}
\end{equation}
We choose $\epsilon_0$ small enough depending only on $C_\sigma$ and $d$, so that the quantity $\sigma(B_Q) - c_Q \mu_{P_Q}(2B_Q)$ is positive and bigger than $\sigma(B_Q)/2 \geq \ell(Q)^{d}/C_\sigma$. Hence, with our choice of $\epsilon_0$, $\alpha_Q$ is bigger than a constant that depends only on $d$ and $C_\sigma$.

By a similar argument, we prove that if $c_Q$ is large and $\dist(\xi_Q,P_Q) \leq 8\ell(Q)$, then $c_Q \mu_{P_Q}(B_Q) - \sigma(2B_Q)$ is positive and bigger than $\ell(Q)^d$. Thus using $\hat f$ in \eqref{a6.2} as before will also implies that $\alpha_Q$ is bigger than a uniform constant.

Assume now that $\dist(2Q,P_{Q}) \leq 5\ell(Q)/\sqrt n$ but $c_Q$ is bigger than $\epsilon_0$. Then by \eqref{prW1}, it means that we can find $x\in P_Q$ such that $\dist(x,\Gamma) \geq 10\ell(Q)/\sqrt n$. We construct the 1-Lipschitz function $f$ on $\R^n$ as
\[ f(X) = \max\{10\ell(Q)/\sqrt n - |X-x|,0\}.\]
The function $f$ is supported in $2^{5}B_Q$, it is a simple consequence of the fact that $Q\subset B_Q$ and $x$ is not far from $Q$. We deduce by definition of $\alpha_\sigma$ (and by our choice of flat measure $\mu_Q$) that 
\[\left| \int f \, d\sigma - \int f \, d\mu_Q \right| \leq 2 [2^{10}\ell(Q)]^{d+1} \alpha_Q.\]
But $f$ do not touch $\Gamma$, so the above estimate becomes
\[c_Q \int_{P_Q} f \, dy \leq 2^{11+10d} \ell(Q)^{d+1}  \alpha_Q, \]
where $dy$ is the $d$-dimensional Lebesgue measure. We can estimate $\int_{P_Q} f\, dy$ with our choice of $f$, we can deduce that 
\begin{equation} \label{muQa}
c_Q \leq C \alpha_Q
\end{equation}
with a constant $C$ that depends only on $d$ and $n$. But we assumed that $c_Q$ is large (bigger than $\epsilon_0$), then \eqref{muQa} implies that $\alpha_Q$ is also bigger than a constant that depends only on $d$ and $C_\sigma$. 
The lemma follows.
\ep

For each $Q \in \cW$, we pick a flat measure. If $\alpha_Q \leq \epsilon$, where $\epsilon$ is the one in Lemma \ref{lemmuQ}, then we take a constant $c_Q$ and a $d$-plane $P_Q$ such that the flat measure $\mu_Q := c_Q \mu_{P_Q}$ satisfies $\dist_Q(\mu_Q,\sigma) \leq 2 \alpha_Q$. If $\alpha_Q \geq \epsilon$, then we take $\mu_Q:= c_Q\mu_{P_Q}$ with $c_Q := 1$ and $P_Q$ any $d$-plane going through $\xi_Q$ and that does not intersect $20Q$ (it is possible since $20Q \not\ni \xi_Q$ is convex). The following properties are proved by Lemma \ref{lemmuQ} when $\alpha_Q$ is small and are immediate by construction when $\alpha_Q$ is large:
\begin{equation} \label{muQ1}
\dist(2Q,P_Q) \geq 5 \ell(Q)/\sqrt n \geq C^{-1} \dist(Q,\Gamma),
\end{equation}
\begin{equation} \label{muQ4}
\dist(\xi_Q,P_Q) \leq 5 \ell(Q),
\end{equation}
\begin{equation} \label{muQ3}
C^{-1} \leq c_Q \leq C,
\end{equation}
and
\begin{equation} \label{muQ2}
\wt \alpha_Q := \dist_Q(\mu_Q,\sigma) \leq C \alpha_Q
\end{equation}
for a constant $C>0$ that depends on $d$, $n$, and $C_\sigma$.

We introduced our choices of flat measure that approximates $\Gamma$, it is now a good time to present the following lemma, that presents how we shall use the rectifiability of $\Gamma$. Since we shall need it later, we set the quantities
\begin{equation} \label{defalphaQk}
\alpha_{Q,k} := \inf_{\mu \in \mathcal F} \dist_{\xi_Q,2^{10+k}\ell(Q)}(\sigma,\mu),
\end{equation}
and we observe that $\alpha_Q$ is $\alpha_{Q,0}$.

\begin{lemma} \label{lemUR}
Let $\Gamma$ be uniformly rectifiable. For $x\in \Gamma$ and $r>0$, define $\cW(x,r)$ as the sub-collection of $\cW$ of cubes $Q$ for which $2Q$ intersects $B(x,r) \subset \R^n$. Then for $x\in \Gamma$, $r>0$, and $k\in \mathbb N$,
\[\sum_{Q \in \cW(x,r)} (\alpha_{Q,k})^2 \ell(Q)^{d} \leq C(C_0+k)r^d,\]
where $C>0$ depends only on $C_\sigma$, $d$ and $n$. By \eqref{muQ2}, we immediately have
\[\sum_{Q \in \cW(x,r)} (\wt \alpha_Q)^2 \ell(Q)^{d} \leq CC_0r^d,\]
where $C$ depends on the same parameters.
\end{lemma}

\bp The proof is pretty much immediate. Let $Q\in \cW(x,r)$. If $y\in \Gamma$ and $s>0$ are such that $|y-\xi_Q| \leq \ell(Q)$ and $2^{11+k}\ell(Q) \leq s \leq 2^{12+k} \ell(Q)$, then $B(y,s) \supset B(\xi_Q,2^{10+k}\ell(Q))$ and the set of functions $Lip(y,s)$ is larger than $Lip(\xi_Q,2^{10+k}\ell(Q))$. The definitions \eqref{a6.2} and \eqref{a6.3} entail that $\alpha_{Q,k} \leq \alpha_\sigma(y,s)$, which can be rewritten
\[(\alpha_{Q,k})^2 \ell(Q)^d \leq C \int_{2^{11+k}\ell(Q)}^{2^{12+k}\ell(Q)} \int_{\Gamma \cap B(\xi_Q,\ell(Q))} |\alpha_\sigma(y,s)|^2 d\sigma(y) \frac{ds}{s}\]
where the constant $C$ depends only on $C_\sigma$. Summing over $Q$ gives that
\begin{equation} \label{URa}
\sum_{Q\in \cW(x,r)} (\alpha_{Q,k})^2 \ell(Q)^d \lesssim \sum_{Q\in \cW(x,r)}\int_{2^{11+k}\ell(Q)}^{2^{12+k}\ell(Q)} \int_{\Gamma \cap B(\xi_Q,\ell(Q))} |\alpha_\sigma(y,s)|^2 d\sigma(y) \frac{ds}{s}
\end{equation}
Notice that the collection 
\begin{equation} \label{URb}
\{(2^{11+k}\ell(Q),2^{12+k}\ell(Q)) \times (\Gamma \cap B(\xi_Q,\ell(Q)))\}_{Q\in \cW}
\end{equation}
is finitely overlapping in $(0,+\infty) \times \Gamma$ (with a uniform constant that depends only on $n$). Indeed, an overlap appears only for cubes $Q,R$ that have the same diameter $D$ (recall that $\cW$ is a collection of dyadic cubes in $\R^n$) and when $|\xi_Q-\xi_R|\leq D$. But the latter implies that $100Q\cap 100R \neq \emptyset$, and given $Q\in \cW$, there is a uniformly finite number of Whitney cubes $R$ that satisfies $\ell(Q) = \ell(R)$ and $100Q\cap 100R \neq \emptyset$.

We use the finite overlap of \eqref{URb} in \eqref{URa} to obtain
\[\sum_{Q\in \cW(x,r)} (\alpha_{Q,k})^2 \ell(Q)^d \lesssim \int_{0}^{(2^{12+k}\sqrt n)r} \int_{B(x,2^{12}\sqrt n r)} |\alpha_\sigma(y,s)|^2 d\sigma(y) \frac{ds}{s}\]
because if a cube $Q$ is in $\cW(x,r)$, then we need to have - for instance - $\diam(Q) \leq \sqrt{n} \, r$ and $\xi_Q\in B(x,100\sqrt n \, r)$. We divide the integral in $s$ into two parts: when $s \leq (2^{12}\sqrt n) r$ and when $(2^{12}\sqrt n) r \leq s \leq (2^{12+k}\sqrt n)r$, and by \eqref{alphabdd} and \eqref{defUR}, we obtain
\[\begin{split}
\sum_{Q\in \cW(x,r)} (\alpha_{Q,k})^2 \ell(Q)^d & \lesssim \int_{0}^{(2^{12}\sqrt n) r} \int_{B(x,2^{12}\sqrt n r)} |\alpha_\sigma(y,s)|^2 d\sigma(y) \frac{ds}{s} \\
& \hspace{3cm} + \int_{(2^{12}\sqrt n) r}^{(2^{12+k}\sqrt n)r} \int_{B(x,2^{12}\sqrt n r)} |\alpha_\sigma(y,s)|^2 d\sigma(y) \frac{ds}{s} \\
& \lesssim \sigma(B(x,2^{12}\sqrt n \, r)) \left[ C_0 + C_\sigma \int_{(2^{12}\sqrt n) r}^{(2^{12}\sqrt n)r}  \frac{ds}{s} \right] \\
& \lesssim r^d [C_0 + k].
\end{split}\]
The lemma follows. \ep

We are almost done with flat measures that approximate $\sigma$. We shall just link a point in $\Omega$ to a flat measure, but we do not need to be gentle, so we take, for $X\in \Omega$, $\mu_X := \mu_Q$ where $Q \in \cW$ is the only dyadic cube containing $X$. Similarly, $c_X$ is $c_Q$ and $P_X$ is $c_Q$ where $X \in Q \in \cW$.
From \eqref{muQ1} and \eqref{muQ3}, it is not very hard to see that we have
\begin{equation} \label{muX1}
\dist(X,P_X) \geq C^{-1} \dist(X,\Gamma).
\end{equation}
and
\begin{equation} \label{muX2}
C^{-1} \leq c_X \leq C
\end{equation}
for a constant $C>0$ that depends only on $d$, $n$, and $C_\sigma$. It is also practical to introduce the alpha numbers relatively to the point $X$ 
\begin{equation} \label{defalphaX}
\alpha(X) := \wt \alpha_Q = \dist_{\xi_Q,2^{10}\ell(Q)}(\sigma,\mu_X),
\end{equation}

The second part of Lemma \ref{lemUR} can now be rewritten as:

\begin{lemma} \label{lemUR2}
If $\Gamma$ is uniformly rectifiable, for every $x\in \Gamma$ and $r>0$, there holds
\begin{equation}\label{UR2b}
\int_{B(x,r)} |\alpha(X)|^2 \dist(X,\Gamma)^{d-n}dX \leq CC_0 r^d,
\end{equation}
where $C>0$ depends only on $C_\sigma$, $d$ and $n$.
\end{lemma}

We are prepared to talk about the soft distance $D_\beta$. We shall also use the vector function defined $\Omega$ by
\begin{equation} \label{defHb}
H_\beta (X) := \int_\Gamma |X-y|^{-d-\beta-1} (X-y) d\sigma(y).
\end{equation}
The purpose of Lemma \ref{Main11a}, with our new notation, is to compare $H_{n-d-1}$ and $D_\beta^{d+1-n} \nabla D_\beta$. 
We will actually prove a more general result, that compares $H_\alpha$ and $D_\beta^{-\alpha} \nabla D_\beta$ for any $\alpha,\beta>0$. Before starting the long computations involving $H_\alpha$ and $D_\beta$, we introduce few notation and make few observations. 

From the definition \eqref{defDb} of $D_\beta$, we can see that the term  $H_\beta(X)$ is immediately bounded by $D_\beta^{-\beta}$ and thus by \eqref{Dbdist}
\begin{equation} \label{bdHb}
|H_\beta| \lesssim D_\beta^{-\beta} \lesssim \dist(X,\Gamma)^{-\beta}.
\end{equation}
Moreover, a direct computation shows that  
\begin{equation} \label{NDbHb}
\nabla D_\beta^{-\beta} = -(d+\beta) H_{\beta +1}
\end{equation}
which entails $\nabla D_\beta = \frac{d+\beta}{\beta} D_\beta^{\beta+1} H_{\beta+1}$  and then
\begin{equation} \label{bdNDb}
|\nabla D_\beta| \lesssim 1.
\end{equation}

Let $c_\beta$ be the number 
\begin{equation} \label{defcb}
c_\beta := \int_{\R^d} (1+|y|^2)^{-\frac{d+\beta}{2}} \, dy.
\end{equation}
Verify that by a simple change of variable that for all $Y,X\in \Omega$
\begin{equation} \label{prcb}
\int |Y-y|^{-d-\beta} \, d\mu_X(y) = c_\beta c_X \dist(Y,P_X)^{-\beta}.
\end{equation}
Let $N_X$ be the unit vector defined for all $X\in \Omega$ by $N_X:=\left[\nabla_Y \dist(Y,P_X)\right]_{Y=X}$.  In one hand, we have by symmetry of $P_X$ that 
\begin{equation} \label{prcb3} \begin{split}
\int |X-y|^{-d-\beta-1}(X-y) \, d\mu_X(y) & = \left(\int |X-y|^{-d-\beta-1}(X-y)\cdot N_X \, d\mu_X(y)\right) N_X \\
& = \left(\int |X-y|^{-d-\beta-1} \dist(X,P_X) \, d\mu_X(y)\right) N_X \\
&= c_{\beta+1} c_X \dist(X,P_X)^{-\beta} N_X,
\end{split}\end{equation}
but also, in the other hand, we can write when $\beta > 1$
\begin{equation} \label{prcb2} \begin{split}
\int |X-y|^{-d-\beta-1}(X-y) \, d\mu_X(y) &= -\frac1{\beta+d-1} \left(\nabla_Y \left[ \int |Y-y|^{-d-\beta+1} \, d\mu_X(y)  \right] \right)_{Y=X} \\ &= \frac{\beta -1}{\beta +d - 1} c_{\beta-1} c_X \dist(X,P_X)^{-\beta} N_X
\end{split}\end{equation}
by \eqref{prcb}.
The combination of the two last estimates gives a nice relation on the coefficients $c_\beta$
\begin{equation} \label{prcb4} 
(\beta + d) c_{\beta+2} = \beta c_\beta \qquad \text{ whenever } \beta >0.
\end{equation}

The next lemma shows the cost of changing the measure $\sigma$ by $\mu_X$ in $D_\beta$ and $H_\beta$.

\begin{lemma} \label{lemDbeta}
Let $\beta>0$. There exists a constant $C>0$ that depends only on $d$, $n$, $\beta$, and $C_\sigma$ such that f or all $X\in \Omega$, 
\begin{equation}\label{Db1}
|D_\beta^{-\beta}(X) - c_\beta c_X \dist(X,P_X)^{-\beta}| \leq C \dist(X,\Gamma)^{-\beta} \left( \alpha(X) +  \sum_{k\in \bN} 2^{-k\beta}\alpha_{Q,k} \right),
\end{equation}
and
\begin{equation}\label{Db1a}
|H_\beta(X) - c_{\beta+1} c_X \dist(X,P_X)^{-\beta} N_X| \leq C \dist(X,\Gamma)^{-\beta} \left( \alpha(X) +  \sum_{k\in \bN} 2^{-k\beta}\alpha_{Q,k} \right),
\end{equation}
where $Q\in \cW$ is the only dyadic cube containing $X$. 
\end{lemma}

\bp The two estimates are proven in the same manner. We shall rigorously prove \eqref{Db1} first, and only explain the differences for \eqref{Db1a}. Let $X\in \Omega$ and let $Q\in \cW$ be the cube containing $X$. We intend to cut the integral
\[D_\beta^{-\beta} = \int_\Gamma |X-y|^{-d-\beta} \, d\sigma(y)\]
into pieces that lives in annuli, so we use cut-off functions that live in the desired annuli. We start by taking a function $\wt \theta_0 :\, \R^n \to \R$ supported in $B(0,2^{10}\ell(Q))$, and satisfying $0\leq \wt \theta_0 \leq 1$ everywhere, $\wt \theta_0 \equiv 1$ on $B(0,2^9\ell(Q))$, and $|\nabla \wt \theta_0| \leq 1/\ell(Q)$. Then we set for $k\geq 1$ the functions
\[\wt \theta_k(y) := \wt \theta_0(2^{-k}y) - \wt \theta_0(2^{-k+1}y),\]
and we translate these functions by taking for all $k\in \bN$
\[\theta_k(y) := \wt \theta_k(y-\xi_Q).\]
The functions $\theta_k$ form a partition of unity, that is
\begin{equation} \label{sumthetak}
\sum_{k \in \bN} \theta_k \equiv 1.
\end{equation}
In addition, for $k\geq 0$
\begin{equation} \label{prtheta1}
\text{ $\theta_k$ is supported in $B_k:= B(\xi_Q,2^{10+k}\ell(Q))$},
\end{equation}
and for $k\geq 1$,
\begin{equation} \label{prtheta2}
\text{ $\theta_k \equiv 0$ on $B_{k-2}$.}
\end{equation}

We shall use the decomposition $D_\beta^{-\beta} = \sum_{k\in \bN} I_k$, where
\begin{equation} \label{1defIk}
I_k := \int_\Gamma |X-y|^{-d-\beta} \theta_k(y) \, d\sigma(y) = \int_\Gamma f_k(y) \, d\sigma(y)
\end{equation}
if $f_k$ is the function defined on $\R^n$ by
\begin{equation} \label{1deffk}
f_k(y) := |X-y|^{-d-\beta} \theta_k(y).
\end{equation}
We intend to approximate $I_k$ by
\begin{equation} \label{1defJk}
J_k := \int_{P_X} f_k \, d\mu_X,
\end{equation}
which can be defined without problems thanks to \eqref{muX1}. The sum of the $J_k$'s can be directly linked to \eqref{Db1}. Indeed, observe that 
\[\sum_{k\in \bN} J_k = \sum_{k\in \bN} \int_{P_X} |X-y|^{-d-\beta} \theta_k(y) \, d\mu_X(y) = \int |X-y|^{-d-\beta} \, d\mu_X(y) = c_\beta c_X \dist(X,P_X)^{-\beta}\]
by \eqref{prcb}. Therefore, we have
\begin{equation} \label{Db2}
|D_\beta^{-\beta}(X) - c_\beta c_X \dist(X,P_X)^{-\beta}| \leq \sum_{k\in \bN} |I_k - J_k|,
\end{equation}
and we just have to estimate the difference $|I_k-J_k|$.

We start with $k=0$. By \eqref{muX1}, $P_X$ stays far from $X$, so we can find $\epsilon\in (0,1)$ independent of $X$ such that $B(X,\epsilon\dist(X,\Gamma))$ intersects neither $P_X$ nor $\Gamma$. If the function $\wt f_0$ is defined as
\begin{equation} \label{Db3}
\wt f_0(y) := \min\{|X-y|^{-d-\beta},[\epsilon\dist(X,\Gamma)]^{-d-\beta}\} \theta_k,
\end{equation}
then $\wt f_0 = f_0$ on $\Gamma \cup P_X$ and we have
\begin{equation} \label{Db4}
I_0 =  \int \wt f_0 \, d\sigma \quad \text{ and } \quad J_0 =  \int \wt f_0 \, d\mu_X.
\end{equation}
By \eqref{prtheta1}, the function $f_0$ is supported in $B_{0} = B(\xi_Q,2^{10}\ell(Q))$, so is our new function $\wt f_0$. But the function $\wt f_0$ is now Lipschitz, with constant $C\dist(X,\Gamma)^{-d-\beta-1}$, or $C\ell(Q)^{-d-1-\beta}$ since $X\in Q$ and $Q$ is a Whitney cube. Hence, 
\begin{equation} \label{Db4a}
|I_0 - J_0| = \left|\int \wt f_0 \, (d\sigma - d\mu_X)\right| \lesssim \ell(Q)^{-\beta} \alpha(X)
\end{equation}
by definition of $\alpha(X)$.

We turn to the case $k\geq 1$. Here, $f_k$ is already Lipschitz. Indeed, it is not hard to check that by construction, we have $X\in B_{-5} = B(\xi_Q, 32\ell(Q))$. So \eqref{prtheta2} entails that $f_k$ is 0 around $X$. Note that 
\begin{equation} \label{Db5}
\text{the Lipschitz constant of $f_k$ is smaller than $C[2^k\ell(Q)]^{-d-\beta-1}$.}
\end{equation}
and
\begin{equation} \label{Db6}
\text{$f_k$ is supported in $B_k = B(\xi_Q,2^{10+k}\ell(Q))$}.
\end{equation}
Those observations will be of use a bit later.

The $\alpha$-numbers will not allow us to compare directly $I_k$ with $J_k$, because $\mu_X$ is not the correct flat measure that approximate $\sigma$ in the ball $B_k$. So for $j\in \bN$, we take flat measures $\mu_{Q,j}=c_{Q,j}\mu_{P_{Q,j}}$ such that 
\[\dist_{\xi_Q,2^{10+j}\ell(Q)}(\sigma,\mu_{Q,j}) \leq 2 \alpha_{Q,j}.\]
The key observation is the fact that 
\begin{equation} \label{Db7}
\dist_{\xi_Q,2^{10+k}\ell(Q)}(\mu_{Q,j-1},\mu_{Q,j}) \lesssim \alpha_{Q,j}.
\end{equation}
for all $1<j\leq k$. The rough idea is that you just need to know the difference between $c_{Q,j}$ and $c_{Q,j-1}$, the distance and the angle between the planes $P_{Q,j}$ and $P_{Q,j-1}$, to be able to bound the Wasserstein distance $\dist_{\xi_Q,2^{10+k}\ell(Q)}(\mu_{Q,j-1},\mu_{Q,j})$. But all these 3 quantities can be estimated with $\alpha_{Q,j}$. A detailed explanation of \eqref{Db7} is given as Appendix \ref{ApWD}. For a similar reason 
\begin{equation} \label{Db8}
\dist_{\xi_Q,2^{10+k}\ell(Q)}(\mu_{Q,0},\mu_{X}) \lesssim \alpha_{Q,0} + \alpha(X).
\end{equation}
We deduce that 
\begin{equation} \label{Db9}
\dist_{\xi_Q,2^{10+k}\ell(Q)}(\mu_{Q,k},\mu_{X}) \lesssim \alpha(X) + \sum_{j=0}^k \alpha_{Q,j}.
\end{equation}

Let us return to the $I_k$ and $J_k$ for $k\geq 1$. We have
\[|I_k - J_k| \leq \left| \int f_k \, d\sigma - \int f_k \, d\mu_{Q,k} \right| + \left| \int f_k \, d\mu_{Q,k} - \int f_k \, d\mu_{X} \right|.\]
The properties \eqref{Db5} and \eqref{Db6} on $f_k$ allows us to bound the above terms with the help of Wasserstein distances. One has
\begin{equation} \label{Db10}\begin{split}
 |I_k - J_k| & \lesssim [2^k\ell(Q)]^{-\beta} \left[ \dist_{\xi_Q,2^{10+k}\ell(Q)}(\sigma,\mu_{Q,k}) + \dist_{\xi_Q,2^{10+k}\ell(Q)}(\mu_{Q,k},\mu_{X}) \right] \\
 & \lesssim [2^k\ell(Q)]^{-\beta} \left[ \alpha(X) + \sum_{j=0}^k \alpha_{Q,j} \right]
 \end{split} \end{equation}
 by the definition of $\mu_{Q,k}$ and by \eqref{Db9}.

The estimates \eqref{Db2}, \eqref{Db4a}, and \eqref{Db10} prove that 
\begin{equation} \label{Db11}\begin{split}
|D_\beta^{-\beta}(X) - c_\beta c_X \dist(X,P_X)^{-\beta}| & \lesssim \ell(Q)^{-\beta} \sum_{k\in \bN} 2^{-\beta k} \left( \alpha(X) + \sum_{j=0}^k \alpha_{Q,j} \right) \\
& \lesssim \ell(Q)^{-\beta} \left( \alpha(X) + \sum_{j\in \bN} 2^{-j\beta} \alpha_{Q,j} \right)
 \end{split} \end{equation}
 by Fubini's theorem. The estimate \eqref{Db1} follows by recalling that, since $X\in Q$ and $Q$ is a Whitney cube, we have $\ell(Q) \eqsim \dist(X,\Gamma)$.
 
 \medskip
 
The estimate \eqref{Db1a} can be established exactly as \eqref{Db1}, by noticing that the argument only requires that the functions inside the integral - namely $|X-y|^{-d-\beta}$ for \eqref{Db1} and $|X-y|^{-d-\beta-2}(X-y)$ for \eqref{Db1a} - are Lipschitz (on $\Gamma \cup P_X$) and have enough decay at infinity.

As before, recall that 
\[H_\beta(X) := \int_\Gamma |X-y|^{-d-\beta-1} (X-y) d\sigma(y)\]
We use the functions 
\[f'_k := |X-y|^{-d-\beta-1} (X-y) \theta_k(y)\]
to make the decomposition $H_\beta = \sum_{k\in \bN} I'_k$ where 
\[I'_k:= \int_\Gamma |X-y|^{-d-\beta-1} (X-y) \theta_k(y) d\sigma(y) = \int_\Gamma f'_k(y) d\sigma(y).\]
We define then $J'_k := \int_{P_X} f'_kd\mu_X$, which satisfies by \eqref{prcb3} that 
\[\sum_k J'_k = c_{\beta+1} c_X \dist(X,P_X)^{-\beta-1} N_X.\] 
And since
\[|H_\beta(X) - c_{\beta+1} c_X \dist(X,P_X)^{-\beta} N_X| \leq \sum_{k\in \bN} |I'_k - J'_k|,\]
we only need to get appropriate bounds on $|I'_k-J'_k|$. We use the same argument as the one given for the proof of \eqref{Db1}, observing that the Lipschitz constant of $f'_k$ is now at most $C[2^k\ell(Q)]^{-d-\beta-1}$ and we obtain that 
\[|I'_k - J'_k| \lesssim [2^k\ell(Q)]^{-\beta} \left[ \alpha(X) + \sum_{j=0}^k \alpha_{Q,j} \right],\]
which implies,
\[|H_\beta(X) - c_{\beta+1} c_X \dist(X,P_X)^{-\beta} N_X| \lesssim  \dist(X,\Gamma)^{-\beta} \left( \alpha(X) + \sum_{j\in \bN} 2^{-j\beta} \alpha_{Q,j} \right)\]
The bound \eqref{Db1a} and the lemma follow.
\ep

For the rest of the section, we take $\alpha,\beta>0$. We overload the notation $\alpha$, that is used both for Tolsa's $\alpha$-number and as a parameter of the quantities $D_\alpha$ and $H_\alpha$, but the two $\alpha$ have such different roles that we believe no confusion should arise from it. But anyway, let us to get rid of any mention to the $\alpha$ as Tolsa number as fast as possible. To that objective, we define the quantity $a(X)$ as
\begin{equation}\label{defoffa}
a(X) := \alpha(X) + \sum_{k\in \bN} 2^{-k\min\{\alpha,\beta\}} \alpha_{Q,k}
\end{equation}
where $Q\ni X$, which is nice because 
\begin{equation} \label{aisCM}
\text{If $\Gamma$ is uniformly rectifiable, then $a$ satisfies the Carleson measure condition.}
\end{equation}
Indeed, by the Cauchy-Schwarz inequality, one has
\[|a(X)|^2 \lesssim |\alpha(X)|^2 + \sum_{k\in \bN} 2^{-k\min\{\alpha,\beta\}} |\alpha_{Q,k}|^2\]
so for every $x\in \Gamma$ and $r>0$, we get
\[\begin{split}
\int_{B(x,r)} |a(X)|^2 \dist(X,\Gamma)^{d-n} dX &\lesssim \int_{B(x,r)} |\alpha(X)|^2 \dist(X,\Gamma)^{d-n} dX \\
& \hspace{3cm} + \sum_{k\in \bN} 2^{-k\min\{\alpha,\beta\}} \sum_{Q\in \cW(x,r)} |\alpha_{Q,k}|^2 \ell(Q)^d \\
& \lesssim C_0r^d (1 + \sum_{k\in \bN} 2^{-k\min\{\alpha,\beta\}} k) \lesssim C_0r^d 
\end{split}\]
where we used the fact that $\dist(X,\Gamma) \eqsim \ell(Q)$ in the first inequality, and Lemmas \ref{lemUR} and \ref{lemUR2} in the second one.
The uniform boundedness of $a$, which is required in the definition of the Carleson measure condition, is a consequence of the fact that $\alpha_\sigma$ - and thus all the quantities constructed from it - is uniformly bounded, as recalled in the introduction.

\medskip

\noindent {\em Proof of Lemma \ref{Main11}.} Assume that $\Gamma$ is uniformly rectifiable. We want to estimate $\nabla [D_\beta/D_\alpha]$, which can be rewritten
\begin{equation} \label{Db13c}
\nabla\left(\frac{D_\beta}{D_\alpha} \right) = \frac{D_\beta}{D_\alpha} \left( \dfrac{\nabla [D_\alpha^{-\alpha}]}{\alpha D_\alpha^{-\alpha}} - \dfrac{\nabla [D_\beta^{-\beta}]}{\beta D_\beta^{-\beta}}  \right) = \frac{D_\beta}{D_\alpha} \left( \dfrac{(d+\beta) H_{\beta+1}}{\beta D_\beta^{-\beta}} - \dfrac{(d+\alpha) H_{\alpha+1}}{\alpha D_\alpha^{-\alpha}}   \right),
\end{equation}
by \eqref{NDbHb}.

Lemma \ref{lemDbeta} infers that
\begin{equation}\label{Db12}
|D_\beta^{-\beta} - c_\beta c_X \dist(X,P_X)^{-\beta}| \leq C \dist(X,\Gamma)^{-\beta} a(X),
\end{equation}
\begin{equation}\label{Db12b}
|H_{\beta+1} -c_{\beta+2} c_X \dist(X,P_X)^{-\beta-1} N_X| \leq C \dist(X,\Gamma)^{-\beta-1} a(X),
\end{equation}
and analogous estimates for  $D_\alpha^{-\alpha}$ and $\nabla[D_\alpha^{-\alpha}]$. So since $D_\alpha \eqsim D_\beta$ by \eqref{Dbdist},
\begin{equation} \label{Db13d}\begin{split}
\left| \nabla\left(\frac{D_\beta}{D_\alpha} \right) \right| \lesssim \left| \dfrac{(d+\beta) H_{\beta +1}}{\beta D_\beta^{-\beta}} - \dfrac{N_X}{\dist(X,P_X)}\right| + \left| \dfrac{(d+\alpha) H_{\alpha+1}}{\alpha D_\alpha^{-\alpha}} - \dfrac{N_X}{\dist(X,P_X)}\right|  
\end{split}\end{equation} 
We shall only bound the first term on the right-hand side above, since the same bound will be obviously obtained on the second term by replacing $\beta$ by $\alpha$. Recall that by \eqref{Dbdist} and \eqref{muX1}
\begin{equation} \label{Db13e}
D_\beta(X) \eqsim \dist(X,\Gamma) \eqsim \dist(X,P_X),
\end{equation} 
the lower bound $\dist(X,P_X) \leq C \dist(X,\Gamma)$ being an easy consequence of \eqref{muQ4}. So we have
\[ \begin{split}
 \left| \dfrac{(d+\beta) H_{\beta +1}}{\beta D_\beta^{-\beta}} -  \dfrac{N_X}{\dist(X,P_X)}\right| 
 & \lesssim \dist(X,P_X)^{\beta} \left| H_{\beta+1}- \dfrac{\beta D_\beta^{-\beta}N_X }{(d+\beta)\dist(X,P_X)}\right| \\
 & \leq \dist(X,P_X)^{\beta} \Big( \left| H_{\beta+1} - c_{\beta+2} c_X N_X \dist(X,P_X)^{-\beta-1}\right|   \\
 & \qquad + \left. \frac{\beta |N_X|}{(d+\beta) \dist(X,P_X)} \left| D_\beta^{-\beta} - \frac{(d+\beta)c_{\beta+2}}{\beta} c_X \dist(X,P_X)^{-\beta}\right|\right)
\end{split}\]
But since $(d+\beta)c_{\beta+2}/\beta = c_\beta$ due to \eqref{prcb4}, we obtain
\[ \begin{split}
 \left| \dfrac{(d+\beta) H_{\beta +1}}{\beta D_\beta^{-\beta}} -  \dfrac{N_X}{\dist(X,P_X)}\right| 
 & \leq \dist(X,P_X)^{\beta} \Big( \left| H_{\beta+1} - c_{\beta+2} c_X N_X \dist(X,P_X)^{-\beta-1}\right|   \\
 & \qquad + \left. \frac{\beta}{(d+\beta) \dist(X,P_X)} \left| D_\beta^{-\beta} - c_\beta c_X \dist(X,P_X)^{-\beta}\right|\right) \\
 & \lesssim \frac{a(X)}{\dist(X,P_X)} \eqsim \dist(X,\Gamma)^{-1} a(X)
\end{split}\]
by \eqref{Db12} and \eqref{Db12b}. Together with \eqref{Db13d}, we conclude that 
\[\dist(X,\Gamma) \left| \nabla\left(\frac{D_\beta}{D_\alpha} \right) \right| \lesssim a(X),\]
which concludes the proof of Lemma \ref{Main11} since $a(X)$ satisfies the Carleson measure condition. \ep

Let us finish the section by the following lemma. Lemma \ref{Main11a} from the introduction is a consequence of the next lemma in the particular case where $\alpha = n-d-1>0$.

\begin{lemma} \label{Main11b}
Let $\Gamma$ be uniformly rectifiable Let $\alpha,\beta >0$.  Then there exist a scalar function $b$ and a vector function $\mathcal V$, both defined on $\Omega$, such that 
\[ H_{\alpha} = (b  \nabla D_\beta + \mathcal V) D_\beta^{-\alpha} \qquad \text{ for } X\in \Omega\]
\begin{enumerate}[(i)]
\item $C_1^{-1} \leq b \leq C_1$,
\item $D_\beta \nabla b \in CM(C_1)$,
\item $|\mathcal V| \leq C_1$,
\item $\mathcal V \in CM(C_1)$,
 \end{enumerate}
 where $C_1$ is a constant that depends only on $C_\sigma$, $C_0$, $\alpha$, $\beta$, $n$, and $d$.
\end{lemma}

\noindent {\em Proof of Lemma \ref{Main11a}.} First, let us explain why we have, for any real $\nu$,
\begin{equation}\label{claimDbcx}
|D_\beta^{\nu}(X) - c_\beta^{-\nu/\beta} c_X^{-\nu/\beta} \dist(X,P_X)^\nu| \leq C \dist(X,\Gamma)^{\nu} a(X).
\end{equation}
By the Mean Value Theorem, if $x,y>0$, we have 
\begin{equation}\label{claimDbcx1}
|x^{-\nu/\beta} - y^{-\nu/\beta}| \leq \left( \sup_{z\in [x,y]} \frac{|\nu|}\beta z^{-\nu/\beta-1}\right) |x-y|.
\end{equation}
We apply the above estimate to $x = D_\beta^{-\beta}/\dist(X,P_X)^{-\beta}$ and $y = c_\beta c_X$. Since $x$ and $y$ are uniformly (in $X$) bounded from above and from below by a positive constant (see \eqref{Db13e}, \eqref{muX2}), we deduce that $\sup_{z\in [x,y]} \frac1\beta z^{-1/\beta-1}$ is bounded uniformly in $X$. As a consequence,
\begin{equation}\label{D-cX<a}
\left| \frac{D_\beta^\nu(X)}{\dist(X,P_X)^\nu} - c_\beta^{-\nu/\beta} c_X^{-\nu/\beta} \right| \lesssim \left| \frac{D_\beta^{-\beta}(X)}{\dist(X,P_X)^{-\beta}} - c_\beta c_X \right| \lesssim a(X)
\end{equation}
thanks to \eqref{Db1}. The claim \eqref{claimDbcx} follows once \eqref{Db13e} is invoked again.

Then we need a smooth substitute for the function $X \to c_X$. We could have been more careful when building $c_X$, so that it would have been smooth from the beginning, but we chose another path. We write $s_X$ for 
\begin{equation} \label{defsx}
s_X := c_1 c_{1/2}^{-2} \frac{D_1(X)}{D_{1/2}(X)},
\end{equation}
where $c_1,c_{1/2}$ are the coefficients defined in \eqref{defcb} for the values $\beta = 1,1/2$, and $D_1,D_{1/2}$ are the smooth distances defined in \eqref{defDb}. We claim that
\begin{equation} \label{cxsx}
|c_X - s_X| \lesssim a(X) \qquad \text{ for } X\in \Omega.
\end{equation}
Indeed, we call \eqref{Db13e} and \eqref{muX2} to justify that 
\begin{multline}
|s_X - c_X|  \lesssim \frac1{\dist(X,\Gamma)} \left| \frac{c_1}{c_X} D_1(X) - c_{1/2}^2 D_{1/2}(X) \right| \\
 \leq  \frac1{\dist(X,\Gamma)} \left( \frac{c_1}{c_X} \left| D_1(X) - c_1^{-1} c_X^{-1} \dist(X,P_X)\right| + c_{1/2}^2 \left| D_{1/2}(X) - c_{1/2}^{-2} c_X^{-2} \dist(X,P_X)\right| \right) \\
 \lesssim a(X)
\end{multline}
by \eqref{claimDbcx}. The claim \eqref{cxsx} follows.

The quantities $s_X$ and $c_X$ are uniformly bounded from above and below by a positive constant, as a consequence of \eqref{Dbdist} and \eqref{muX2}. So if we mimic the proof of \eqref{claimDbcx}, using the Mean Value Theorem, we get that for any real number $\nu$, one  has
\begin{equation} \label{cxsx2}
|c_X^{\nu} - s_X^\nu| \lesssim a(X).
\end{equation}
We can now replace the $c_X$ by $s_X$ in the estimates \eqref{claimDbcx} and \eqref{Db2}. Indeed, 
\[\begin{split}
|D_\beta^{\nu}(X) - c_\beta^{-\nu/\beta} s_X^{-\nu/\beta} \dist(X,P_X)^\nu| 
 \leq |D_\beta^{\nu}(X)  - c_\beta^{-\nu/\beta} & c_X^{-\nu/\beta} \dist(X,P_X)^\nu| \\
& + c_\beta^{-\nu/\beta} \dist(X,P_X)^\nu | c_X^{-\nu/\beta} - s_X^{-\nu/\beta}|
\end{split}\]
So by \eqref{claimDbcx}, \eqref{Db13e}, and \eqref{cxsx2}, for any power $\nu \in \R$ and any $\beta >0$.
\begin{equation} \label{Dbsx}
|D_\beta^{\nu}(X) - c_\beta^{-\nu/\beta} s_X^{-\nu/\beta} \dist(X,P_X)^\nu| \lesssim \dist(X,\Gamma)^\nu a(X).
\end{equation}
Similarly, we have for $\beta >0$,
\begin{equation} \label{Hbsx}
|H_\beta(X) - c_{\beta+1} s_X \dist(X,P_X)^{-\beta} N_X| \lesssim \dist(X,\Gamma)^{-\beta} a(X).
\end{equation}

\medskip

We finished our preliminary estimates. We set 
\[b(X) :=  \frac{\beta c_{\alpha+1}}{(d+\beta)c_{\beta+2}} [c_\beta s_X]^{(\beta + 1 - \alpha)/\beta}.\]
Since $\nabla b$ is - up to a constant - $s_X^{(1-\alpha)/\beta} \nabla s_X$, and since $s_X$ is the quotient $D_1/D_{1/2}$, Lemma \ref{Main11} and \eqref{Dbdist} give the conclusions $(i)$ and $(ii)$ of the lemma under proof. So it remains to check that 
\[\mathcal V:= D_\beta^{\alpha} H_\alpha - b\nabla D_\beta\]
is uniformly bounded and satisfies the Carleson measure condition.  The uniform bound on $\mathcal V$ is easy, and is a consequence of the fact that $|H_\alpha| \lesssim D_\alpha^{-\alpha} \eqsim D_\beta^{-\alpha}$ and $|\nabla D_\beta| \lesssim 1$ (see \eqref{bdHb} and \eqref{bdNDb}).

We are left with the proof of $(iv)$, i.e. the fact that $\mathcal V \in CM$. Thanks to \eqref{NDbHb}, we have that
\begin{equation} \label{estV1} \begin{split}
|\mathcal V| & = \left| D_\beta^{\alpha} H_\alpha - \frac {d+\beta}{\beta} b D_\beta^{\beta+1} H_{\beta+1} \right| \\
& \leq D_\beta^{\beta+1} \Big( \left| D_\beta^{\alpha-\beta-1} H_\alpha - c_{\alpha+1} [c_\beta]^{(\beta+1-\alpha)/\beta} [s_X]^{1 + \frac{\beta + 1 - \alpha}{\beta}} \dist(X,P_X)^{-\beta-1} N_X  \right|  \\
& \qquad + \left. \left| \frac{d+\beta}\beta b H_{\beta+1} - c_{\alpha+1} [c_\beta]^{(\beta+1-\alpha)/\beta} [s_X]^{1 + \frac{\beta + 1 - \alpha}{\beta}} \dist(X,P_X)^{-\beta-1} N_X \right| \right) \\
& := D_{\beta}^{\beta+1} (V_1 + V_2).
\end{split}\end{equation}
We start by bounding $V_2$. We use the expression of $b$ to get
\begin{equation} \label{estV2} \begin{split}
V_2 & = \frac{c_{\alpha+1}}{c_{\beta+2}} [c_\beta s_X]^{(\beta+1-\alpha)/\beta}  |H_{\beta+1} - c_{\beta+2} s_X \dist(X,P_X)^{-\beta-1} N_X | \\
& \lesssim \dist(X,\Gamma)^{-\beta-1} a(X) \lesssim D_\beta^{-\beta-1} a(X).
\end{split}\end{equation}
by \eqref{muX2}, \eqref{Hbsx}, and then \eqref{Dbdist}. As for $V_1$, we write
\begin{equation} \label{estV3} \begin{split}
V_1 & = D_\beta^{\alpha-\beta-1} \left|H_\alpha - c_{\alpha+1} [c_\beta]^{(\beta+1-\alpha)/\beta} [s_X]^{1 + \frac{\beta + 1 - \alpha}{\beta}} \dist(X,P_X)^{-\beta-1} D_\beta^{\beta+1-\alpha} N_X \right| \\
& \leq D_\beta^{\alpha-\beta-1} \Big(\left| H_\alpha - c_{\alpha+1} s_X \dist(X,P_X)^{-\alpha} N_X \right| \\
& \qquad + c_{\alpha+1} [c_\beta]^{(\beta+1-\alpha)/\beta} [s_X]^{1 + \frac{\beta + 1 - \alpha}{\beta}} \dist(X,P_X)^{-\beta-1} \\
&
\hspace{5cm} \left| D_\beta^{\beta+1-\alpha} - [c_{\beta}s_X]^{(\alpha - \beta-1)/\beta} \dist(X,P_X)^{\beta+1-\alpha} \right| \Big) \\
& \lesssim D_\beta^{\alpha-\beta-1} (\dist(X,\Gamma)^{-\alpha} + \dist(X,P_X)^{-\beta-1} \dist(X,\Gamma)^{\beta+1-\alpha} ) a(X) \\
& \lesssim D_\beta^{-\beta-1} a(X)
\end{split}\end{equation}
by \eqref{Dbsx}, \eqref{Hbsx}, and then \eqref{Db13e}. We combine the estimates \eqref{estV1}--\eqref{estV3} to deduce that 
\[|\mathcal V| \lesssim a(X).\]
Since by \eqref{aisCM}, the quantity $a$ satisfies the Carleson measure condition, $|\mathcal V|$ also satisfies the the Carleson measure condition. Conclusion $(iv)$ and the lemma follow.
\ep

\appendix

\section{Wasserstein distances between planes} \label{ApWD}

We shall use the notation introduced in the beginning of Subsection \ref{SS1.2}. In particular, $\mathcal F$ is the set of the flat measures $\mu$ that can be written $\mu := c \mu_{P}$, where $c>0$ is a constant and $\mu_{P}$ is the ($d$-dimensional) Lebesgue measure on an affine plane $P \subset \R^n$ of dimension $d$.

Given two flat measures $\mu_1 = c_1 \mu_{P_1}$ and $\mu_2 = c_2 \mu_{P_2}$ whose supports intersect $B(x_0,r/2)$, we want to estimate
\begin{equation} \label{distP1P2}
\dist_{x_0,r}(\mu_1 ,\mu_2) := r^{-d-1} \sup_{f \in Lip(x_0,r)} \left| \int f \, d\mu_1 - \int f \, d\mu_2 \right|.
\end{equation}
The let the reader check that by definition of $Lip(x_0,r)$, we necessarily have 
\begin{equation} \label{Ap1}
\dist_{x_0,r}(\mu_1 ,\mu_2) \lesssim c_1 + c_2.
\end{equation}

The plane $P_2$ may be orthogonal to $P_1$, which in our mindset means that the orthogonal projection of $P_2$ onto $P_1$ is a strict subset $P_{11}$ of $P_1$. In this case, we can find a point $y\in B(x_0,3r/4) \cap P_1$ for which $\dist(y,P_{11}) > r/4$ and hence $\dist(y,P_2)>r/4$. Using the function $f(x) = \max\{ r/4 - |y-x|, 0\}$ to estimate $\dist_{x_0,r}(\mu_1 ,\mu_2)$, we deduce that $\dist_{x_0,r}(\mu_1,\mu_2) \gtrsim c_1$. By symmetry of the role of $\mu_1$ and $\mu_2$, we also deduce that $\dist_{x_0,r}(\mu_1,\mu_2) \gtrsim c_2$. Altogether, we have
\begin{equation} \label{Ap2}
\dist_{x_0,r}(\mu_1 ,\mu_2) \eqsim c_1 + c_2,
\end{equation}
which is close to the worst case scenario by \eqref{Ap1}.

Without loss of generality, we can also assume that $c_2 \leq c_1$. By translation and rotation invariance, we can also assume that $x_0\in B(0,r/2)$ and that
$$P_1 = \R^d = \{(x,0) \in \R^n, \, x\in \R^d\}.$$
 If $P_1$ and $P_2$ are not orthogonal, then $P_2$ can be written as
 $$P_2 := \{(x,xa+b) \in \R^n, \, x\in \R^d\}$$
 where $a$ is a $d\times (n-d)$ matrix and $b\in \R^{n-d}$.

\begin{lemma} \label{lemAp1}
$$\dist_{x_0,r}(\mu_1,\mu_2) \lesssim c_1 (|a| + \frac{|b|}{r}) + (c_1 - c_2).$$
\end{lemma}

In the above lemma, $|a|$ denotes a matrix norm, the exact definition has no importance since all matrix norms are equivalent.

\medskip

\bp
We write $c_2(a)$ for the quantity
$$c_2(a) = c_2 \left(1 + \sum_{i,j} a_{i,j}^2 \right)^{-1/2}.$$
The integral over $\mu_2$ in \eqref{distP1P2} can be rewritten as
$$\int f \, d\mu_2 = c_2(a) \int_{\R^d} f(x,xa+b) dx.$$
Therefore,
\[\begin{split}
\left| \int f \, d\mu_1 - \int f \, d\mu_2 \right| 
& \leq c_1 \left| \int [f(x,0) - f(x,xa+b)] dx \right| + |c_1 - c_2(a)| \left| \int f(x,xa+b)  dx \right|.
\end{split}\]
We use the fact that $f$ is $1$-Lipschitz and supported in $B(x,r)$ - in particular we have $f(y) \leq r$ for all $y \in \R^n$ - to obtain that 
\[\begin{split}
\left| \int f \, d\mu_1 - \int f \, d\mu_2 \right| 
& \lesssim c_1 (|a|r +|b|)r^d + (c_1 - c_2(a)) r^{d+1} \\
& \lesssim c_1 (|a|r +|b|)r^d + (c_1 - c_2) r^{d+1} + (c_2 - c_2(a)) r^{d+1}.
\end{split}\]
By definition of $c_2(a)$, we have that $c_2 - c_2(a) \lesssim c_2 |a| \lesssim c_1 |a|$. So we deduce that 
$$\left| \int f \, d\mu_1 - \int f \, d\mu_2 \right|  \lesssim c_1 (|a|r +|b|)r^d + (c_1 - c_2) r^{d+1}.$$
The lemma follows from taking the supremum over all the function $f\in Lip(x,r)$, and then by dividing by $r^{d+1}$.
\ep

\begin{lemma}  \label{lemAp2}
If $|a| \leq 1$, then 
$$\dist_{x_0,r}(\mu_1,\mu_2) \eqsim c_1 (|a| + \frac{|b|}{r}) + (c_1 - c_2).$$
If $|a| \geq 1$, then
$$\dist_{x_0,r}(\mu_1,\mu_2) \eqsim c_1.$$
\end{lemma}

\bp
We start with the case $|a|\geq 1$. In this case, it is roughly the same as the case where $P_1$ and $P_2$ are orthogonal. Wince $a$ is large enough, we can find a value $y \in P_1 \cap B(0,r/4) \subset B(x,3r/4)$ such that $\dist(y,P_2) \geq \alpha r$ for some $\alpha \in (0,1)$ that depends only on the matrix norm we used when we wrote $|a|$.

We take then the function $f(x) = \max\{0, \alpha r - |y-x| \}$ in \eqref{distP1P2} and we deduce that 
$$\dist_{x_0,r}(\mu_1 ,\mu_2) \geq c_1.$$
Since we assumed that $c_2 \leq c_1$, \eqref{Ap1} gives us the reverse inequality.

\medskip

We consider now that $|a| \leq 1$. We proved the upper bound in Lemma \ref{lemAp1}, so we only need to check the lower bound. We prove first that 
\begin{equation} \label{Ap3}
\dist_{x_0,r}(\mu_1 ,\mu_2) \geq c_1 - c_2.
\end{equation}
Indeed, construct the function $\phi : \R \to \R$ such that $\phi \equiv 1$ on $[-3/2,3/2]$, $\phi$ is supported in $(-2,2)$, and $\phi$ is 1-Lipschitz. We build then $f(x,y)$ on $\R^d \times \R^{n-d}$ as
$$f(x,y) = \frac{r}{16} \phi\left(\frac{8|y|}{r}\right) \phi\left(\frac{k|x|}{r}\right)$$
with a constant $k \leq 8$ to find. Observe that we already have that $f$ is supported in $B(0,r/2) \subset B(x_0,r)$ and is $1$-Lipschitz, so in particular $f\in Lip(x_0,r)$. Since $|a|\leq 1$, we can take $k \gtrsim 1$ such that $f(x,y) = \frac r{16}\phi(k|x|/r)$ whenever $(x,y) \in P_2$. With this function $f$, we have
\[\begin{split}
\left| \int f \, d\mu_1 - \int f \, d\mu_2 \right| & = \left| c_1 \int_{\Rd} f(x,0) dx - c_2(a) \int_{\R^d} f(x,xa+b) dx  \right| \\
& = \int_{\R^d} r \phi(k|x|/r) (c_1-c_2(a)) dx \\
& \gtrsim r^{d+1} (c_1 - c_2(a)) \geq r^{d+1} (c_1-c_2).
\end{split}\]
In the above computations, $c_2(a)$ denotes the constant defined in the beginning of the proof of Lemma \ref{lemAp1}. The claim \ref{Ap3} follows.

If $\frac{10|b|}{r} \geq |a|$, then we consider the function 
$$f(x,y) = \frac{|b|}{200} \phi\left(\frac{100|y|}{|b|}\right) \phi\left(\frac{100|x|}{r}\right),$$
which belongs to $Lip(x,r)$ and takes the value 0 on $P_2$. Using this function in \eqref{distP1P2}, we deduce that 
\begin{equation} \label{Ap4}
\dist_{x_0,r}(\mu_1 ,\mu_2) \gtrsim c_1|b|/r  \eqsim c_1 (|a| + \frac{|b|}{r}).
\end{equation}

If $\frac{10|b|}{r} \leq |a|$, then we can find a point $z\in P_1 \cap B(0,r/4)$ such that $|za+b| \geq |a|r/10$. We use then the function
$$f(x,y) =  \frac{|a|r}{200} \phi\left(\frac{100|y|}{|a|r}\right) \phi\left(\frac{100|x-z|}{r}\right).$$
Observe that $f \in Lip(x,r)$ and $f \equiv 0$ on $P_2$. Consequently,
\begin{equation} \label{Ap5}
\dist_{x_0,r}(\mu_1 ,\mu_2) \gtrsim c_1 |a|   \eqsim c_1 (|a| + \frac{|b|}{r}).
\end{equation}

The combination of \eqref{Ap3}, \eqref{Ap4}, and \eqref{Ap5} leads to the first part of the lemma.
\ep

\begin{remark}
Lemma \ref{lemAp2} entails that the function $r \mapsto \dist_{x_0,r}(\mu_1,\mu_2)$ is essentially decreasing as long as $r$ is large enough. More precisely, 
\[ \dist_{x_0,2^kr} (\mu_1,\mu_2) \leq C \dist_{x_0,r} (\mu_1,\mu_2) \]
if both $P_1$ and $P_2$ intersect $B(x_0,r/2)$. 

A consequence of this result is \eqref{Db7}, that we used in the proof of Lemma \ref{lemDbeta}.
\end{remark}


\begin{thebibliography}{AAA}

\bibitem[Azz]{Azzam} J. Azzam. {\em Semi-uniform domains and the $A_\infty$ property for the harmonic measure.} Int. Math. Res. Not. (2021), no. 9, 6717--6771.

\bibitem[AHMNT]{AHMNT} J. Azzam, S. Hofmann, J.M.  Martell, K. Nystr{\"o}m,  T. Toro. 
{\em A new characterization of chord-arc domains}, JEMS {\bf 19} (2017), no. 4, 967--981.

\bibitem[AHM3TV]{AHM3TV} J. Azzam, S. Hoffman, M. Mourgoglou,  J. M. Martell, S. Mayboroda, X.  Tolsa, A. Volberg. 
{\em Rectifiability of harmonic measure}, Geom. Funct. Anal., {\bf 26} (2016), no. 3, 703--728.

\bibitem[AHMMT]{AHMMT} J. Azzam, S. Hofmann, J. M. Martell, M. Mourgoglou, X. Tolsa. {\em Harmonic measure and quantitative connectivity: geometric characterization of the $L^p$ solvability of the Dirichlet problem.} Invent. Math. {\bf 222} (2020), no. 3, 881--993.

\bibitem[Ba]{Ba} M.\,Badger. {\it 
Null sets of harmonic measure on NTA domains: Lipschitz approximation revisited.}
Math. Z. 270 (2012), no. 1-2, 241--262. 

\bibitem[BJ]{BJ} C. Bishop, P. Jones.
{\it Harmonic measure and arclength,}
Ann. of Math. (2), {\bf 132} (1990), 511--547.

\bibitem[CFK]{CFK} L. Caffarelli, E. Fabes, C. Kenig. {\em Completely singular elliptic-harmonic measures.} Indiana Univ. Math. J. {\bf 30} (1981), no. 6, 917--924.

  \bibitem[Da]{Da}
B.~E.~J. Dahlberg, \emph{Estimates of harmonic measure}.  Arch. Rational Mech.
  Anal. \textbf{65} (1977), no.~3, 275--288. 

\bibitem[DJ]{DJ} G. David, D. Jerison.  
{\it Lipschitz approximation to hypersurfaces, harmonic measure, and singular integrals.} 
Indiana Univ. Math. J., {\bf 39} (1990), no. 3, 831--845.

\bibitem[DS1]{DS1} G. David, S. Semmes.
{\em Singular integrals and rectifiable sets in $\re^n$: Beyond Lipschitz graphs.}
Asterisque, {\bf 193} (1991).

\bibitem[DS2]{DS2} G. David, S. Semmes.
{\em Analysis of and on uniformly rectifiable sets.} 
Mathematical Surveys and Monographs, {\bf 38}. American Mathematical Society, Providence, RI, 1993.

\bibitem[DEM]{DEM} G. David, M. Engelstein, S. Mayboroda. {\em Square functions, non-tangential limits and harmonic measure in co-dimensions larger than 1}. Duke Math. J. {\bf 170} (2021), no. 3, 455--501.

\bibitem[DFM1]{DFMprelim} G. David, J.  Feneuil, S. Mayboroda. 
{\em Elliptic theory for sets with higher co-dimensional boundaries.} Mem. Amer. Math. Soc. {\bf 274} (2021), no. 1346, vi+123 pp.

\bibitem[DFM2]{DFMAinfty} G. David, J.  Feneuil, S. Mayboroda. 
{\em Dahlberg's theorem in higher co-dimension.} J. Funct. Anal. 276 (2019), no. 9, 2731--2820.

\bibitem[DFM3]{DFMKenig} G. David, J.  Feneuil, S. Mayboroda. 
{\em A new elliptic measure on lower dimensional sets.} Acta Math. Sin. (Engl. Ser.) 35 (2019), no. 6, 876--902.

\bibitem[DFM4]{DFMprelim2} G. David, J.  Feneuil, S. Mayboroda. {\em Elliptic theory in domains with boundaries of mixed dimension}. Preprint, arxiv:2003.09037.

\bibitem[DM]{DM} G. David, S. Mayboroda. {\em Harmonic measure is absolutely continuous with respect to the Hausdorff measure on all low-dimensional uniformly rectifiable sets}. Preprint, arXiv:2006.14661.

\bibitem[DPP]{DPP2015} M. Dindo\v s, S.\,Petermichl, J.\,Pipher, {\it BMO solvability and the $A_\infty$ condition for second order parabolic operators}. Ann. Henri Poincar\'e (C) Non Linear Analysis {\bf 34} (2017), no. 5, 1155--1180.

\bibitem[GT]{GT} D. Gilbarg, N. S. Trudinger.
{\em Elliptic partial differential equations of second order.}
Classics in Mathematics. Springer-Verlag, Berlin (2001). Reprint of the 1998 edition.

\bibitem[HLMN]{HLMN} S. Hofmann, P. Le, J. M. Martell, K. Nystr\"om. {\it The weak-$A_\infty$ property of harmonic and $p$-harmonic measures implies uniform rectifiability}. Anal. PDE 10 (2017), no. 3, 513--558.

\bibitem[HLM]{HLM} S. Hofmann, P. Le, and A. Morris, 
{\em Carleson measure estimates and the Dirichlet problem for degenerate elliptic equations}. Anal. PDE 12 (2019), no. 8, 2095--2146.

\bibitem[HM]{HM1}
S. Hofmann, J.M. Martell. 
{\em Uniform rectifiability and harmonic measure {I}: uniform rectifiability implies {P}oisson kernels in {$L^p$}}. Ann. Sci. \'Ec. Norm. Sup\'er. (4), {\bf 47} (2014), no. 3, 577--654.
  
\bibitem[HMU]{HMU} S. Hofmann, J.M. Martell, I. Uriarte-Tuero. 
{\em Uniform rectifiability and harmonic measure, {II}: {P}oisson kernels in {$L^p$} imply uniform rectifiability.} 
Duke Math. J., {\bf 163} (2014), no. 8, 1601--1654.

\bibitem[JK]{JK} D. Jerison, C. Kenig. {\em The Dirichlet problem in nonsmooth domains.} Ann. of Math. (2), {\bf 113} (1981), no. 2, 367--382.

\bibitem[Ken]{KenigB} C. E. Kenig.
{\em Harmonic analysis techniques for second order elliptic boundary value problems.} 
CBMS Regional Conference Series in Mathematics, {\bf 83}. Amer. Math. Soc., Providence, RI, 1994.

\bibitem[KKiPT]{KKiPT} C. Kenig, B. Kirchheim, J. Pipher, T. Toro. 
{\it Square Functions and the $A_\infty$ Property of Elliptic Measures.} 
J. Geom. Anal., {\bf 26} (2016), no. 3, 2383--2410.

\bibitem[KP]{KePiDrift} C. Kenig, J. Pipher.  
{\it The Dirichlet problem for elliptic equations with drift terms.} 
Publ. Mat., {\bf 45} (2001), no. 1, 199--217. 

\bibitem[LN]{LN} J.\,Lewis, K.\,Nystr\"om. {\it Quasi-linear PDEs and low-dimensional sets,} JEMS {\bf 20} (2018), no. 7, 1689--1746. 

\bibitem[Lv]{Lv} M. Lavrent'ev. 
{\em Boundary problems in the theory of univalent functions.} 
Amer. Math. Soc. Transl. (2), {\bf 32} (1963), 1--35.

\bibitem[MZ]{MZ} S. Mayboroda, Z. Zihui. {\em Square function estimates, BMO Dirichlet problem, and absolute continuity of harmonic measure on lower dimensional sets.} Anal. PDE 12 (2019), 1843--1890.

\bibitem[MM]{MM} L. Modica, S. Mortola. {\em Construction of a singular elliptic-harmonic measure.} Manuscripta Math. {\bf 33} (1980/81), no. 1, 81--98.
  
 \bibitem[Pi]{PipherICM} J. Pipher. {\em Carleson Measures and elliptic boundary value problems.} Proceedings of the ICM 2014.

\bibitem[RR]{RR} F. \& M. Riesz.
{\it \"Uber die randwerte einer analtischen funktion.}
Compte Rendues du Quatri\`eme Congr\`es des Math\'ematiciens Scandinaves, Stockholm 1916,
Almqvists and Wilksels, Upsala, 1920.

\bibitem[Russ]{Russ} E. Russ. {\em The atomic decomposition for tent spaces on spaces of homogeneous type}. CMA/AMSI Research Symposium ``Asymptotic Geometric Analysis, Harmonic Analysis, and Related Topics'', Proc. Centre Math. Appl. Austral. Nat. Univ. {\bf 42} (2007), 125--135.

\bibitem[Se]{Se} S. Semmes. 
{\it Analysis vs. geometry on a class of rectifiable hypersurfaces in $\R^{n}$.} 
Indiana Univ. Math. J., {\bf 39} (1990), no. 4, 1005--1035.

\bibitem[Ste]{Stein93} E. M. Stein.
{\em Harmonic analysis: real-variable methods, orthogonality, and oscillatory integrals}.
Princeton Mathematical Series, {\bf 43}. Princeton University Press, Princeton, N.J., 1993.

\bibitem[Tol]{Tolsa09} X. Tolsa.
{\em Uniform rectifiability, {C}alder\'on-{Z}ygmund operators with odd kernel, and quasiorthogonality.}
{Proc. Lond. Math. Soc. (3)}, {\bf 98} (2009), no. 2, 393--426.

\bibitem[Tor]{ToroICM} T. Toro. {\em Potential Analysis meets Geometric Measure Theory}. Proceedings of the ICM 2010, III.

\bibitem[Wu]{Wu} J.-M. Wu. {\it On singularity of harmonic measure in space.} Pacific J. Math., {\bf 121} (1986), no. 2, 485--496.

\bibitem[Z]{Z} W.P. Ziemer.
{\it Some remarks on harmonic measure in space.} 
Pacific J. Math., {\bf 55} (1974), 629--637.



\end{thebibliography}
\end{document}